\documentclass[12pt]{article}
\usepackage[english]{babel}
\usepackage[latin1]{inputenc}
\usepackage{amsfonts,amssymb,amsmath, epsfig}
\usepackage{color,graphicx,graphics,psfrag}
\usepackage{amsmath,amstext,amssymb,amsfonts, amscd}
\usepackage{amsthm}
\usepackage{vmargin} 
\usepackage[only,llbracket,rrbracket]{stmaryrd}

\usepackage{hyperref} 

\newcommand{\veps}{\varepsilon}

\def\Box{\leavevmode\vbox{\hrule
     \hbox{\vrule\kern4pt\vbox{\kern4pt}%
           \vrule}\hrule}}
\def\blackbox{\leavevmode\vrule height 5pt width 4pt depth 0pt\relax}
\def\endproof{\null\hfill {$\blackbox$}\bigskip}

\def\paragraph#1{{\bf #1\ }}

\newtheorem{lemma}{Lemma}[section]  

\newtheorem{theorem}[lemma]{Theorem}

\newtheorem{definition}[lemma]{Definition}

\newtheorem{proposition}[lemma]{Proposition}

\newtheorem{remark}{Remark}[section]

\title{Phase transition and diffusion among socially interacting self-propelled agents}      

\begin{document}
\maketitle

\centerline{\scshape Alethea B.T. Barbaro }
\medskip
{\footnotesize
 \centerline{Department of Mathematics}
   \centerline{Case Western Reserve University}
   \centerline{ 10900 Euclic Avenue--Yost Hall Room 220}   
   \centerline{ Cleveland, OH 44106-7058, USA.}
   \centerline{alethea.barbaro@case.edu}
} 

\medskip

\centerline{\scshape Pierre Degond}
\medskip
{\footnotesize
 \centerline{ Universit\'e de Toulouse; UPS, INSA, UT1, UTM}
   \centerline{Institut de Math\'ematiques de Toulouse}   
   \centerline{31062 Toulouse, France.}
      \centerline{and}
   \centerline{CNRS; Institut de Math\'ematiques de Toulouse UMR 5219}
      \centerline{31062 Toulouse, France.}
       \centerline{pierre.degond@math.univ-toulouse.fr}
}

\bigskip

\begin{abstract}
We consider a hydrodynamic model of swarming behavior derived from the kinetic description of a particle system combining a noisy Cucker-Smale consensus force and self-propulsion. In the large self-propulsion force limit, we provide evidence of a phase transition from disordered to ordered motion which manifests itself as a change of type of the limit model (from hyperbolic to diffusive) at the crossing of a critical noise intensity. In the hyperbolic regime, the resulting model, referred to as the `Self-Organized Hydrodynamics (SOH)', consists of a system of compressible Euler equations with a speed constraint. We show that the range of SOH models obtained by this limit is restricted. To waive this restriction, we compute the Navier-Stokes diffusive corrections to the hydrodynamic model. Adding these diffusive corrections, the limit of a large propulsion force yields unrestricted SOH models and offers an alternative to the derivation of the SOH using kinetic models with speed constraints. 
\end{abstract}

\medskip
\noindent
{\bf Acknowledgements:} This work has been supported by the french `Agence Nationale pour la Recherche (ANR)' 
in the frame of the contracts `MOTIMO' (ANR-11-MONU-009-01) and `CBDif-Fr' (ANR-08-BLAN-0333-01).

\medskip
\noindent
\textbf{Key words:} Swarm, Cucker-Smale model, Vicsek model, self-propulsion, hydrodynamic model, diffusion, Chapman-Enskog expansion.

\medskip
\noindent
\textbf{AMS subject classification:} 35L60, 35K55, 35Q80, 82C05, 82C22, 82C70, 92D50.

\section{Introduction}
\label{sec:into}

There is a considerable literature devoted to the observation and understanding of systems of swarming agents. Examples of such systems in nature are fish schools \cite{Becco_etal_PhysicaA06, Katz_etal_PNAS11}, bird flocks \cite{Ballerini_etal_PNAS08, Lukeman_etal_PNAS10}, insect swarms \cite{Buhl_etal_Science2006, Couzin_Franks_PRoySocB03} or migrating cell assemblies \cite{Yamao_etal_PLOSone11} (see also the reviews \cite{Couzin_Krause_03, Vicsek_Zafeiris}). Some simple inanimate physical systems also exhibit collective behavior \cite{Deseigne_etal_PRL10, Suematsu_etal_PRE10}. Many of the models proposed in the literature are `Individual-Based Models (IBM)'. They consist in following the dynamics of the agents and their interactions over the course of time. The `three zone model' of Aoki \cite{Aoki, Couzin_etal_JTB02} postulates that interactions obey long-range attraction, short-range repulsion and medium-range alignment. Vicsek et al. \cite{Vicsek_etal_PRL95} have proposed a simplified version of this model where particles move at a constant speed and interact through alignment only. In spite of its simplicity, the Vicsek model exhibits complex features which have triggered a large literature (see \cite{Vicsek_Zafeiris} for a review and references). On the other hand, the Cucker-Smale model \cite{Cucker_Smale_IEEE07} is based on a large-scale velocity consensus formation and does not impose any constraint on the particle speed. The Cucker-Smale model has triggered considerable mathematical activity \cite{Carrillo_etal_SIMA10, Cucker_Mordecki_JMPA08, Ha_etal_CMS09, Ha_Liu_CMS09, Ha_Tadmor_KRM08, Shen_SIAP07}. Many other kinds of IBM's of collective motion can be found and it is impossible to cite them all (see e.g. \cite{Carrillo_etal_M3AS10, Chuang_etal_PhysicaD07, Dorsogna_etal_PRL06, Mogilner_etal_JMB03} and the review \cite{Vicsek_Zafeiris}). Additionally, comparisons of models with data can be found e.g. in \cite{Barbaro_etal_MCS09, Barbaro_etal_ICESJMS09, Hemelrijk_Hildenbrandt_PLOSone11}.

IBM's are very successful but become computationally intensive for large systems. For this reason, macroscopic models of fluid type have been proposed in the literature. Macroscopic models of collective motion have been derived from heuristic rules and symmetry considerations in \cite{Toner_Tu_PRL95, Toner_etal_AnnPhys05}. The rigorous derivation of continuum models usually starts from a statistical version of the IBM, the so-called kinetic model. Kinetic models of collective motion have been derived in \cite{Bolley_etal_M3AS11, Bolley_etal_AML12, Carrillo_etal_SIMA10} for various versions of the Cucker-Smale and Vicsek models. The convergence of the kinetic Cucker-Smale model to the kinetic Vicsek model is shown in \cite{Bostan_Carrillo}. In \cite{Bertin_etal_PRE06, Bertin_etal_JPA09}, a Boltzmann type kinetic model has been proposed for binary collision processes which mimics the Vicsek alignment dynamics. In the same spirit, a Boltzmann-Povzner type approach which mimics the Cucker-Smale process and its fluid limit has been developed in \cite{Fornasier_etal_PhysicaD11}. In \cite{Bertin_etal_PRE06, Bertin_etal_JPA09}, a hydrodynamic model for the binary version of the Vicsek interaction is derived from the kinetic model under an assumption of weak anisotropy of the velocity distribution function. In \cite{Ratushnaya_etal_PhysicaA07} a direct passage from the Vicsek IBM to a fluid model is attempted. The first derivation of hydrodynamic-like equations from the mean-field kinetic version of the Vicsek model has been performed in \cite{Degond_Motsch_M3AS08}. Further elaboration of the model can be found in \cite{Degond_etal_submitted1, Degond_etal_submitted2, Frouvelle_M3AS12, Frouvelle_Liu_SIMA12}. Diffusive corrections to the model of \cite{Degond_Motsch_M3AS08} have been derived in \cite{Degond_Yang_M3AS10} and bear analogies with the model proposed in \cite{Mishra_etal_PRE10}. Other kinds of macroscopic models can be found in \cite{Bertozzi_etal_Nonlinearity09, Eftimie_etal_PNAS07, Mogilner_EK_JMB99, Topaz_Bertozzi_SIAP04, Topaz_etal_BMB06}.

The aim of the present work is twofold. The first objective is to give evidence of a phase transition from disordered to ordered motion in a hydrodynamic model of socially interacting agents with self-propulsion. Specifically, we want to emphasize the role of the self-propulsion in the emergence of the phase transition. Such evidence has been given for the first time in \cite{Toner_Tu_PRL95}. However, the model of \cite{Toner_Tu_PRL95} is based on analogy with the Vicsek IBM, not on an actual derivation from it. In \cite{Toner_Tu_PRL95}, the techniques used to show the emergence of a phase transition are rather complex: they are based on stability analysis in the linear case and renormalization group theory in the nonlinear case. In the present work, the model we investigate is derived from a simple IBM of collective motion combining a noisy Cucker-Smale consensus force and self-propulsion. The phase transition appears simply when analyzing the behavior of the model in the limit of a large self-propulsion force. It manifests itself as a change in the type of limit model at the crossing of a critical noise intensity. Above the critical noise, the limit model is of diffusion type, while below the critical noise, it is of hyperbolic type. To the best of our knowledge, the present work is the first instance where phase transitions in particle swarms have been evidenced in this way from a hydrodynamic model. A similar approach, but at the level of the kinetic model, can be found in  \cite{Degond_etal_submitted1, Frouvelle_Liu_SIMA12}. 

The second goal of this work is to discuss the relative merits of the Cucker-Smale and the Vicsek models for the derivation of hydrodynamic models of particle swarms. As mentioned above, the Vicsek kinetic model imposes that kinetic velocities be of constant norm, while no such constraint exists in the Cucker-Smale model. Instead, a self-propulsion force is imposed to force the particle speed to stay close to a `comfort' velocity. In \cite{Bostan_Carrillo}, it is shown that the Cucker-Smale model relaxes to the Vicsek model when the intensity of the self-propulsion force tends to infinity. The derivation of hydrodynamic models from the Vicsek kinetic model is considerably complexified by the velocity norm constraint. Indeed, momentum conservation is lost and the use of conserved quantities (or collision invariants), which is the cornerstone of the derivation of hydrodynamic models, cannot be implemented. In \cite{Degond_Motsch_M3AS08}, this problem has been overcome by the introduction of a new concept of `Generalized Collision Invariant'. But if the hydrodynamic limit could be performed equivalently on the Cucker-Smale model, these unpleasant technicalities would be proven unnecessary. 

Unfortunately, performing the hydrodynamic limit on the Cucker-Smale model and then letting the self-propulsion force tend to infinity is not equivalent to performing the hydrodynamic limit on the Vicsek model. This is the second main result of the present work. Indeed, the type of the resulting model is the same, but the coefficients of the model are not the same. Specifically, the limit model has the form of a system of isothermal compressible Euler equations for the swarm density $\rho$ and the mean velocity direction $\omega$ (also referred to as the polarization field, see e.g. \cite{Mishra_etal_PRE10}). The velocity direction $\omega$ is a vector of norm one. To maintain this constraint, the model includes some non-conservative terms. Additionally, the convective derivatives involved in the mass and momentum transport are {\bf not the same}, a signature of a loss of Galilean invariance (see e.g. \cite{Tu_etal_PRL98}). This model has been referred to in \cite{Degond_etal_submitted2} as the `Self-Organized Hydrodynamic (SOH)' model. In the present work, we show that the SOH model derived from the Cucker-Smale kinetic equation necessarily involves {\bf the same} convective derivatives in the density and in the velocity equations. Therefore, with the Cucker-Smale model, we cannot access the whole range of possible hydrodynamic limits that we can access with the Vicsek model. With the Cucker-Smale model, we only get a sub-class of these models, which limits its practical applications: having different convective derivatives for $\rho$ and $\omega$ increases the likelihood of correctly reproducing emergent phenomena in swarms, such as cluster formation, waves, etc.

This restriction, which is a significant disadvantage of the Cucker-Smale approach, can be weakened, at least partially. Indeed, in the last part of the present work, we include small diffusive corrections to the hydrodynamic limit of the Cucker-Smale model, by means of a Chapman-Enskog method (see e.g. \cite{Degond_Birkhauser04} for a review). If the self-propulsion is taken to infinity in the resulting compressible Navier-Stokes system, then a more general SOH model, in particular with possibly different convective derivatives for $\rho$ and $\omega$, can be derived. This approach is limited by the necessity to keep the diffusive corrections small; this limits the range of the coefficients of the SOH model which can be obtained. This paper shows that most SOH models can be realized in a fairly general context as hydrodynamic limits of either Cucker-Smale or Vicsek kinetic models.  For these reasons, we conclude that the two approaches are somehow equivalent in the amount of technical work: while the Vicsek model made use of generalized collision invariants, starting from the Cucker-Smale model necessitates dealing with the complex diffusion terms.

The paper is organized as follows. The problem is set up in Section \ref{sec:setting}. Then, some functional properties of the operators are recalled in Section \ref{S:properties_Q}. The hydrodynamic and large self-propulsion limits are derived in Section \ref{sec:hydro}. The diffusive corrections are dealt with in Section \ref{sec:diffusive}. Finally, a conclusion and some perspectives are drawn in Section \ref{sec_conclu}. Three appendices collect some of the more technical proofs.

\section{Setting of the problem}
\label{sec:setting}

\subsection{Velocity consensus in self-propelled agent systems}
\label{subsec_ibm}

We consider a system of agents with positions $x_i(t) \in {\mathbb R}^d$ and velocity $v_i(t) \in {\mathbb R}^d$, where $d$ is the system dimension (in practice equal to $2$ or $3$),  $t \geq 0$ is the time and $i \in \{1, \ldots, N\}$ is the agent's label. These agents are subject to a self-propulsion force which tends to restore a comfort velocity $a>0$, and to a social force which drives them to the average velocity of the neighboring agents. Addtionally, they are subject to random velocity fluctuations which account for potential misperceptions and their propensity to leave the swarm and explore a new environment. The equations of motion are given by 
\begin{eqnarray}
& & \hspace{-1cm} \dot x_i = v_i, \quad  \quad dv_i = {\mathcal F}_i \, dt + \tau^{-1} \big( 1 - \frac{|v_i|^2}{a^2} \big) v_i \, dt + \sqrt{2D} \, dB_t^i, \label{eq:ibmxv}
\end{eqnarray}
with 
\begin{eqnarray}
& & \hspace{-1cm} {\mathcal F}_i = \sigma^{-1} (\bar v_i - v_i), \quad  \quad \bar v_i = \frac{\sum_{j=1}^N K(|x_j-x_i|) v_j}{\sum_{j=1}^N K(|x_j-x_i|)} .  \label{eq:ibmF}
\end{eqnarray}
The social force ${\mathcal F}_i $ is written in (\ref{eq:ibmF}) as a relaxation force towards the average velocity $\bar v_i$ in the neighborhood of particle $i$. The relaxation rate is $\sigma^{-1}$ (in other words, $\sigma$ is the typical time needed for agent $i$ to align with the velocity of his neighbors). The kernel $K$, supposed spherically symmetric for the sake of simplicity, describes how the various partner velocities $v_j$ are combined according to the distance of $j$ to $i$. For instance, if $K$ is the indicator function of the ball of radius $R$, it means that the agents adopt the mean velocity of the other agents within a distance $R$. The second term in the expression of $dv_i$ in (\ref{eq:ibmxv}) is the self-propulsion force. It takes the form of a relaxation term driving $|v_i|$ towards $a$ at rate $\tau^{-1}$. In other words, it takes a time $\tau$ for the velocity $|v_i|$ to relax to the comfort speed $a$. Finally, the last term is the velocity fluctuation term, where $B_t^i$ are independent normalized Brownian processes and $D>0$ is the diffusion coefficient. The force ${\mathcal F}_i $ has been previously proposed by Cucker and Smale \cite{Cucker_Smale_IEEE07} as a model for consensus formation in particle swarms. A noisy version of the Cucker-Smale model is proposed in \cite{Cucker_Mordecki_JMPA08}.

In the large particle limit $N \to \infty$, by adapting the arguments in \cite{Bolley_etal_M3AS11}, the empirical measure of the system
$$ \mu_t^N (x,v) = \frac{1}{N} \sum_{i=1}^N \delta_{(x_i(t),v_i(t))}(x,v), $$
where $\delta_{(x_i(t),v_i(t))}(x,v)$ is the Dirac delta at $(x_i(t),v_i(t))$, can be approximated by a continuous distribution function $f(x,v,t)$. It solves the following Fokker-Planck equation: 
\begin{equation}
\partial_t f + \nabla_x \cdot (v f) = -\nabla_v \cdot \left[ \mathcal{F}_{f}f + \tau^{-1} \big( 1-\frac{|v|^2}{a^2} \big)vf \right] + D \Delta_v f, \label{eq:FPdim}
\end{equation}
with
\begin{eqnarray}
&& \hspace{-1cm} \mathcal{F}_f = \sigma^{-1} (\bar{v}_f - v), \quad \quad 
\bar{v}_f = \frac{\int K({|x-y|})\omega f(y,\omega,t)\,dy \,d\omega}{\int K({|x-y|}) f(y,\omega,t)\,dy \,d\omega } . 
\label{eq:Fdim}
\end{eqnarray}
The left-hand side expresses particle displacement at velocity $v$. The right-hand side consists of three terms. The first one is the consensus force. The second term is the self-propulsion force. The last term takes into account the random velocity fluctuations.

\subsection{Scaling}
\label{subsec_scaling}

In order to understand the roles of the various terms, it is useful to introduce dimensionless quantities. We set $x_0$ and $t_0$ to be space and time units and deduce units of velocity $v_0 = x_0/t_0$ and force ${\mathcal F}_0= x_0/t_0^2$. We assume that the range of the interaction kernel $K$ is $R$, meaning that we can write $ K(|x|) = \tilde K(|x|/R)$ with $\tilde K$ having second moment of order $1$ (i.e. $\int \tilde K(|\tilde x|) |\tilde x|^2 \, d\tilde x = O(1)$; we assume that $\tilde K$ is normalized to $1$, i.e. $\int \tilde K(|\tilde x|) \, d\tilde x = 1$). We now introduce dimensionless variables $\tilde x = x/x_0$, $\tilde t = t/t_0$, $\tilde v = v/v_0$ and the change of variables $\tilde f (\tilde x, \tilde v, \tilde t) = x_0^d \, v_0^d \, f(x_0 \tilde x, v_0 \tilde v, t_0 \tilde t)$. Finally, we introduce the dimensionless parameters: 
\begin{eqnarray}
&& \hspace{-1cm} \hat R = \frac{R}{x_0}, \quad \quad \hat \sigma = \frac{\sigma}{t_0}, \quad \quad \hat a = \frac{a}{v_0}, \quad \quad \hat \tau = \frac{\tau}{t_0}, \quad \quad \hat D = \frac{t_0}{v_0^2} D. \label{eq:dimpar}
\end{eqnarray}
In this new system of coordinates, the system is written: 
\begin{eqnarray}
&& \hspace{-1cm} \partial_t f + \nabla_x \cdot (v f) = -\nabla_v \cdot \left[ \mathcal{F}_{f}f + \hat \tau^{-1} \big( 1-\frac{|v|^2}{\hat a^2} \big)vf \right] + \hat D \Delta_v f, \label{eq:FPadim} \\
&& \hspace{-1cm} \mathcal{F}_f = \hat \sigma^{-1} (\bar{v}_f - v), \quad \quad 
\bar{v}_f = \frac{\int K\big(\frac{|x-y|}{\hat R} \big) \, \omega \,  f(y,\omega,t)\,dy \,d\omega}{\int K\big(\frac{|x-y|}{\hat R} \big) \, f(y,\omega,t)\,dy \,d\omega } ,  \label{eq:Fadim}
\end{eqnarray}
where we have dropped the tildes for the sake of clarity. Now, by fixing the relations between the five dimensionless parameters (\ref{eq:dimpar}), we define the regime we are interested in. We suppose that the diffusion and social forces are simultaneously large, while the range of the social force tends to zero. The parameters of the self-propulsion are kept of order 1. More specifically, we let $\varepsilon \ll 1$ be a small parameter and we assume that $\hat D = {\mathcal O}(1/\varepsilon)$ (large diffusion), $\hat \sigma^{-1} = {\mathcal O}(1/\varepsilon)$ (large social force), $\hat R = {\mathcal O}(\varepsilon)$ (small range of social interaction), while $\hat \tau^{-1} = {\mathcal O}(1)$ and $\hat a = {\mathcal O}(1)$ (parameters of the social force are order unity). In order to highlight these scaling assumptions, we define constants $D^\sharp$, $\sigma^\sharp$, $R^\sharp$, which are all ${\mathcal O}(1)$ and such that 
\begin{eqnarray}
&& \hspace{-1cm} \hat D = \frac{1}{\varepsilon} D^\sharp, \quad \quad \hat \sigma = \varepsilon \sigma^\sharp, \quad \quad \hat R = \varepsilon R^\sharp. \label{eq:parscaled}
\end{eqnarray}
Then, with these new notations, and dropping all `hats' and `sharps', we get the following scaled system: 
\begin{eqnarray}
&&\hspace{-1cm} \veps \left[ \partial_t f^\veps + \nabla_x \cdot (v f^\veps) + \tau^{-1} \nabla_v \cdot \Big( \big(1-\frac{|v|^2}{a^2} \big)v f^\veps \Big) \right]=  -\nabla_v \cdot \left( \mathcal{F}^\veps_{f^\veps}f^\veps \right) + D \Delta_v f^\veps,
\label{eq:rescaledkin} \\
&&\hspace{-1cm} {\mathcal F}^\veps_{f} = \sigma^{-1} (\bar{v}^\veps_{f} - v), \quad 
\bar{v}_f^\veps = \frac{\int K(\frac{|x-y|}{\veps R})\omega f(y,\omega,t)\,dy \,d\omega}{\int K(\frac{|x-y|}{\veps R}) f(y,\omega,t)\,dy \,d\omega}. \label{vfeps}
\end{eqnarray}
We will investigate the limit as  $\varepsilon\to 0$ of this system, while all other parameters (i.e. $\tau$, $a$, $D$, $\sigma$ and $R$) are kept fixed. Hence, we highlight the dependence of $f$ upon $\varepsilon$.  

We can simplify the problem by using Taylor's expansion. At leading order, we find:
\begin{align}
\bar{v}_f^\veps &=  u_f + O(\veps^2) , \quad \quad
u_f = \frac{\int f \omega \,d\omega}{\int f \,d\omega}. \label{vfeps2}
\end{align}
Then, we have
$$ \mathcal{F}^\veps_{f^\veps}  =  \sigma^{-1} (u_{f^\veps} - v) + O(\veps^2) , $$
which leads to:
\begin{eqnarray} 
& & \hspace{-1cm} \veps \left[ \partial_t f^\veps + \nabla_x \cdot (v f^\veps) + \tau^{-1} \nabla_v \cdot \Big( \big(1-\frac{|v|^2}{a^2} \big)v f^\veps \Big) \right]= \nonumber \\
& & \hspace{2cm} = -\nabla_v \cdot \left( \sigma^{-1} (u_{f^\veps} - v) f^\veps  \right) + D \Delta_v f^\veps + O(\veps^2) .
\label{E:ourPDE}
\end{eqnarray}
We will drop the $O(\varepsilon^2)$ terms, since, as a first step,  we consider only the leading and first order terms. Then, problem (\ref{E:ourPDE}) can be written as 
\begin{align} \label{E:ourPDE_1}
\veps \left[ \partial_t f^\veps + \nabla_x \cdot (v f^\veps) + \tau^{-1} \nabla_v \cdot \Big( \big(1-\frac{|v|^2}{a^2} \big)v f^\veps \Big) \right] = Q(f^\veps), 
\end{align}
with the collision operator $Q(f)$:
\begin{equation} \label{E:collisionAlphaZero}
Q(f) = \nabla_v \cdot \left[ \sigma^{-1} (v-u_f) f   + D \nabla_v f \right]. 
\end{equation}

Some remarks concerning scaling (\ref{eq:parscaled}) can be made. The diffusion and social forces are supposed of the same order of magnitude and much larger than all other forces. They counterbalance each other. Indeed, the social force makes the agents adopt the same velocity while diffusion tends to spread the velocities out. This balance results in a Maxwellian velocity profile (i.e. Gaussian in velocity space), as shown later on. The choice which is made here is to assume that the self-propulsion force is weaker. Another choice would have been to make the self-propulsion force as large as the social force and the diffusion. In this case, the balance would involve three different effects and would result in more complicated equilibria. This investigation is in progress \cite{Barbaro_etal_inprep}. The interaction range is supposed to tend to zero like the inverse of the interaction rate. It is no surprise that, at leading order, only spatially local interaction terms remain (which can be seen in (\ref{E:collisionAlphaZero}) by the replacement of the non-local average velocity $\bar{v}_f^\veps$ by the local mean velocity $u_f$). Again, other choices can be made. In \cite{Degond_etal_submitted2}, in the case of the Vicsek model (which is the limit of (\ref{E:ourPDE_1}) when $\tau \to 0$), it is shown that the different choice $\hat R = {\mathcal O}(\sqrt \varepsilon)$ leads to a different macroscopic limit when $\varepsilon \to 0$. This choice takes better care of the non-local character of the interaction. It will be investigated in future work.

Our plan is now to investigate the hydrodynamic limit $\veps \to 0$ in this model. To this end, we first examine the properties of the collision operator $Q$.

\section{Properties of $Q$} 
\label{S:properties_Q}

When $\varepsilon \to 0$ in (\ref{E:ourPDE_1}), $f^\varepsilon$ formally converges to an element of the null-space of $Q$, i.e. a function $f$ such that $Q(f) = 0$. In this limit, the dynamics are characterized by the projection of the left-hand side of (\ref{E:ourPDE_1}) onto the space orthogonal to the range of $Q$. This space is spanned by the so-called `collision invariants'. In this section, we successively determine the null-space of $Q$ and its collision invariants.

\subsection{Null-Space}
\label{SS:null}

We first define the Maxwellian with mean velocity $u \in {\mathbb R}^d$ and temperature $T = \sigma D>0$ as follows:
\begin{equation}
M_u (v):= \frac{1}{(2 \pi T)^\frac{d}{2}}\exp \Big(- \frac{|v-u|^2}{2 T}\Big).
\label{eq:norm_Maxw}
\end{equation}  
Note that $M_u$ satisfies $\int M_u \,dv = 1$ and $\int M_u v \,dv = u$.  

To proceed, we need to determine a functional setting.  Let $u \in {\mathbb R}^d$ and define a weighted $L^2$-space $H_u$ such that
\begin{equation*}
H_u := \{ \phi : \int \phi^2 \, M_u dv < + \infty \} 
\end{equation*}
and a weighted $H^1$-space $V_u$ such that
\begin{equation*}
V_u := \{ \phi \in H_u: \int |\nabla_v \phi|^2 M_u dv < + \infty\}
\end{equation*}
with the associated norms:
\begin{align}
&|\phi|^2_{H_u} := \int {\phi^2}{M_u} \, dv, \quad |\phi|^2_{V_u} := \int |\nabla_v \phi|^2 \, M_u \, dv, \quad ||\phi||_{V_u}^2 = |\phi|_{H_u}^2 + |\phi|_{V_u}^2, \label{eq:inner}
\end{align}
and inner products $(\cdot, \cdot)_{H_u}$, $(\cdot, \cdot)_{V_u}$, and $((\cdot, \cdot))_{V_u}$ respectively.

\begin{lemma}
(i)The operator $Q$ given by \eqref{E:collisionAlphaZero} can be reformulated as:
\begin{equation} \label{E:collision2}
Q(f) = D \nabla_v \cdot \left[ M_{u_f} \nabla_v \big( \frac{f}{M_{u_f}} \big) \right].
\end{equation}

\noindent (ii) For any function $f(v)$ such that $f/M_{u_f} \in V_{u_f}$ and for any function $g \in V_{u_f}$, we have:
\begin{align}
\int Q(f) \, g \,dv = - D \int M_{u_f} \nabla_v \big(\frac{f}{M_{u_f}} \big) \cdot \nabla_v g \,dv. \label{eq:self0}
\end{align}
In particular, we have: 
\begin{align}
\int Q(f)\big( \frac{f}{M_{u_f}} \big) \,dv = -D \int M_{u_f} \left|\nabla_v \big(\frac{f}{M_{u_f}} \big) \right|^2 \,dv. \label{eq:self}
\end{align}

\noindent (iii) The null-space Ker$ \, Q = \{ f(v) \; | \; f/M_{u_f} \in V_{u_f} \mbox{ and } Q(f) = 0 \}$ is given by:
\begin{equation}
\text{Ker} \; Q = \{ \rho M_u \; | \;  \rho \geq 0, \, \,  u \in \mathbb{R}^d \}. 
\label{eq:KerQ}
\end{equation}

\label{lem:Q}
\end{lemma}

\medskip
\noindent
The proof of this lemma is postponed to Appendix \ref{app_prop_Q}. An element $\rho M_u$ of Ker$(Q)$ is called a local thermodynamic equilibrium with density $\rho$ and mean velocity $u$.

\subsection{Collision invariants}
\label{subsec:CI}

\begin{definition}
A function $\psi(v)$ is said to be a \emph{collision invariant} (CI) if and only if 
$$ \int Q(f) \psi \,dv = 0, $$ 
for every $f$ such that $f/M_{u_f} \in V_{u_f}$ and $\psi \in V_{u_f}$. The set of CI's is denoted by ${\mathcal C}$. It is a vector space. 
\label{def:CI}
\end{definition}

\medskip
\noindent
We have the: 

\begin{proposition}
We have $ {\mathcal C} = \text{Span}\{ 1, v_1, \dots, v_d \}$.
In other words,  $\psi$ is a CI if and only if there exists  $a\in \mathbb{R}$ and $b \in \mathbb{R}^d$ such that $\psi(v)  = a+b \cdot v$.
\label{prop:CI}
\end{proposition}

\medskip
\noindent
The proof of this proposition is again postponed to Appendix \ref{app_prop_Q}.

\setcounter{equation}{0}
\section{Hydrodynamic limit and fast relaxation}
\label{sec:hydro}

The goal of this section is to investigate the formal limit $\varepsilon \to 0$ in (\ref{E:ourPDE_1}) and to examine some of the properties of the limit system relative to the propulsion force. More precisely, we exhibit a phase transition when the intensity of the velocity fluctuations crosses a certain threshold dependent on the magnitude of the propulsion velocity.

\subsection{Hydrodynamic limit}
\label{subsec_hydro}
 
The goal of this section is to prove the following formal theorem: 

\begin{theorem}
Let $f^\veps$ be the solution of equation (\ref{E:ourPDE_1}) associated to an initial datum $f_I(x,v)$. We suppose that $f_I$ is independent of $\varepsilon$ for simplicity. We assume that solutions of  (\ref{E:ourPDE_1}) exist on any time interval $[0,{\mathcal T}]$. Assume that $f^\veps \to f^0$ as $\veps \to 0$ as smoothly as needed, which means in particular that derivatives of $f^\veps$ converge to the corresponding derivatives of $f^0$. Then, there exist two functions $\rho(x,t) >0$ and $u(x,t) \in {\mathbb R}^d$ such that 
\begin{equation}
f^0 (x,v,t)= \rho(x,t) M_{u(x,t)}(v), \quad \forall (x,v,t) \in {\mathbb R}^{2d} \times [0,{\mathcal T}]. 
\label{eq:f0}
\end{equation}
Furthermore, $\rho$ and $u$ satisfy the following system of isothermal compressible Euler equations with relaxation: 
\begin{align}
&\partial_t \rho + \nabla_x \cdot (\rho u ) = 0 , \label{E:alpha1}\\
&\partial_t (\rho u ) + \nabla_x \cdot (\rho u \otimes u ) + T \nabla_x \rho = -\frac{1}{\tau} \rho u \big( \frac{|u|^2 + (d+2) T}{a^2} - 1 \big) , \label{E:alpha2}
\end{align}
associated with initial data $(\rho_I, u_I)$ such that 
$$ \rho_I = \int f_I \, dv, \quad \rho_I u_I= \int f_I \, v \, dv. $$
\label{thm:hydro}
\end{theorem}

\medskip
\noindent
{\bf Proof.} First, we note from (\ref{E:ourPDE_1}) that $Q(f^\varepsilon) = O(\varepsilon)$. Therefore $Q(f^0) = 0$, which, because of (\ref{eq:KerQ}), implies that $f^0$ is of the form (\ref{eq:f0}).

Next, using the CI's given by Proposition \ref{prop:CI}, we multiply (\ref{E:ourPDE_1}) successively by $1$ and $v$ and use the fact that the right-hand side vanishes upon integration. Then, we get the following conservation relations, which are valid for any $\varepsilon$: 
\begin{eqnarray} \label{E:massConservation}
& & \partial_t \rho^\veps + \nabla_x \cdot j^\veps = 0, \\
\label{E:conservationOfMomentum}
& & \partial_t j^\veps + \nabla_x \cdot \Sigma^\veps = \tau^{-1} q^\varepsilon , 
\end{eqnarray}
with $\rho^\varepsilon$, $j^\varepsilon$, $\Sigma^\varepsilon$, the density, flux and pressure tensor associated to $f^\varepsilon$, given by: 
\begin{equation} \rho^\veps = \int f^\varepsilon \, dv, \quad j^\veps = \int f^\varepsilon \, v \, dv , \quad \Sigma^\veps = \int f^\varepsilon \, (v \otimes v) \, dv, 
\label{eq:rhojSeps}
\end{equation}
and right-hand side $q^\varepsilon$ given by 
\begin{equation} q^\veps = \int f^\varepsilon v (1 - \frac{|v|^2}{a^2}) \, dv. 
\label{eq:qeps}
\end{equation}

Now, letting $\varepsilon \to 0$, we can express $j = \lim j^\varepsilon$, $\Sigma = \lim \Sigma^\varepsilon$ and $q = \lim q^\varepsilon$ as functions of $\rho$ and $u$: 
\begin{eqnarray*}
&&  j = \int \rho M_u \, v \, dv = \rho u ,\\
&& \Sigma = \int \rho M_u \, (v \otimes v) \, dv = \rho (u \otimes u) + \rho T \mbox{Id}, \\
&& q = \int \rho M_u \, v (1 - \frac{|v|^2}{a^2}) \, dv = \rho u (\frac{|u|^2 + (d+2) T}{a^2} -1) .\\
\end{eqnarray*}
Inserting these expressions into the conservation equations leads to (\ref{E:alpha1}), (\ref{E:alpha2}). The statement about the initial conditions is obvious. \endproof

\begin{remark}
Dividing \eqref{E:alpha2} by $\rho$ and using \eqref{E:alpha1}, the momentum conservation equation can be written in non-conservative form:
\begin{equation} \label{E:nonconservativeForm}
(\partial_t + u \cdot \nabla_x ) u +  T \frac{\nabla_x \rho}{\rho} = - \frac{1}{\tau a^2} u \big(|u|^2 - (a^2 -  (d+2) T) \big).
\end{equation}
\label{rem:nonconservative}
\end{remark}

\noindent
This model is nothing but the Euler system of isothermal compressible gas dynamics with a forcing term. The mass conservation equation (\ref{E:alpha1}) (also known as the continuity equation) has a standard form. The momentum balance equation (\ref{E:alpha2}) involves an isothermal pressure $T \rho$ on the left-hand side and a self-propulsion force on the right-hand side. The temperature $T = \sigma D$ is proportional to the ratio of the intensities of the velocity fluctuations $D$ and of the social force $\sigma^{-1}$. The self-propulsion force takes the form of a relaxation of the fluid velocity to a comfort fluid velocity $\sqrt{a^2 -  (d+2) T}$. In comparison to the force acting on individual particles, the force acting on the fluid involves a term depending on the temperature. The temperature being a truly macroscopic quantity, this term cannot have any counterpart at the particle level. Here, at the fluid level, the comfort fluid velocity becomes a pure imaginary number when $T$ is larger than a critical temperature $T_c = a^2/(d+2)$. When the temperature crosses $T_c$, a phase transition occurs. This phase transition is studied in the next section when the intensity $\tau^{-1}$ of the self-propulsion force is taken to infinity, a limit which we refer to as the `fast relaxation limit'.

\subsection{Fast relaxation limit in the hydrodynamic model}
\label{subsec:tauto0}

We recall that, according to our notation, $\tau$ is the ratio of the physical relaxation time to the time unit $t_0$ (see Section \ref{subsec_scaling}) and is a dimensionless parameter. The goal of this section is to investigate the limit $\tau \to 0$ in the hydrodynamic model (\ref{E:alpha1}), (\ref{E:alpha2}).  We denote by $(\rho^\tau, u^\tau)$ the solution for finite $\tau$ and assume that its limit $(\rho,u)$ as $\tau \to 0$ exists and is as smooth as needed. We will use the non-conservative form of the model, which we recall here for the sake of convenience:
\begin{align}
&\partial_t \rho^\tau + \nabla_x \cdot (\rho^\tau u^\tau ) = 0 , \label{E:alpha1_del}\\
&(\partial_t + u^\tau \cdot \nabla_x ) u^\tau +  T \nabla_x \ln \rho^\tau = - \frac{1}{\tau} u^\tau (\frac{|u^\tau|^2 + (d+2) T}{a^2} -1).  \label{E:alpha2_del}
\end{align}
Letting $\tau \to 0$ formally in (\ref{E:alpha2_del}) leads to $|u|^2 + (d+2)T - a^2 = 0$. Therefore, there are two cases according to whether the quantity $(d+2)T - a^2$ is positive or negative, i.e. according to the position of $T$ with respect to the critical temperature $T_c$ defined by
\begin{equation}
T_c = \frac{a^2}{d+2}. 
\label{eq:Tth}
\end{equation}

\subsubsection{Case $T>T_c$ (large noise)}
\label{subsubsec:large_noise}

We let $(d+2)T - a^2 = (d+2) (T - T_c) := s^2 >0$. The constant $s$ only depends on the problem data and not on the solution. In this case, equation  \eqref{E:alpha2_del} 
becomes
\begin{equation}
(\partial_t + u^\tau \cdot \nabla_x) u^\tau + T \nabla_x \ln \rho^\tau = - \frac{1}{\tau a^2} u^\tau (|u^\tau|^2 + s^2) .
\label{eq:utau_strong}
\end{equation}
To examine the limit $\tau \to 0$, we need to find the equilibria of the right-hand side of (\ref{eq:utau_strong}), i.e. the solutions $u$ of $u(|u|^2 + s^2) = 0$. Obviously, the only solution is $u=0$. Additionally, at least in the spatially homogeneous setting, this is a stable solution. Indeed, the unique solution of 
$$ \frac{du}{dt} = - u (|u|^2 + s^2), \quad u(0) = u_0 ,$$
satisfies $u(t) \to 0$ as $t \to \infty$. So, in the spatially non-homogeneous case, we formally have $u^\tau \to 0$ as $\tau \to 0$. Therefore, the formal limit of (\ref{E:alpha1}), (\ref{E:alpha2}) gives 
$$ \partial_t \rho = 0, \quad u=0.$$ 
In order to get a more precise description of the limit $\tau \to 0$, we need to rescale time and velocity. 

With this aim, we let $t' = \tau t$ and $u^\tau (x,t) = \tau \tilde u^\tau(x,t')$, $\rho^\tau (x,t) = \tilde \rho^\tau (x,t')$.  Inserting this into (\ref{E:alpha1}), (\ref{E:alpha2}), we find (dropping the tildes):
\begin{align}
&\partial_t \rho^\tau + \nabla_x \cdot (\rho^\tau u^\tau ) = 0 , \label{E:alpha1_til}\\
& \tau^2 (\partial_t + u^\tau \cdot \nabla_x) u^\tau + T \nabla_x \ln \rho^\tau = -\frac{1}{a^2 } u^\tau (\tau^2 |u^\tau|^2 + s^2) . \label{E:alpha2_til}
\end{align}
The behavior of this system when $\tau \to 0$ is that of a diffusion. More precisely, we state:

\begin{proposition}
Assume that the solution $(\rho^\tau, u^\tau)$ of system (\ref{E:alpha1_til}), (\ref{E:alpha2_til}) is smooth and converges smoothly towards a pair $(\rho,u)$ as $\tau \to 0$. Then, $\rho$ satisfies the following diffusion equation:
\begin{equation}
\partial_t \rho - D_{\mbox{\scriptsize diff}} \Delta_x \rho = 0, \quad \quad D_{\mbox{\scriptsize diff}} = \frac{TT_c}{T-T_c}. 
\label{eq:heat}
\end{equation}
and $u = - D_{\mbox{\scriptsize diff}} \nabla_x \ln \rho$. 
\label{prop_large_noise}
\end{proposition}

\medskip
\noindent
{\bf Proof.}
Taking $\tau \to 0$ in equation (\ref{E:alpha2_til}), we get $u = -\frac{a^2 T}{s^2} \nabla_x \ln \rho = - \frac{TT_c}{T-T_c} \nabla_x \ln \rho$.  The limit model therefore follows from the continuity equation (\ref{E:alpha1_til}) and is given by (\ref{eq:heat}). \endproof

\begin{remark}
We note that $D_{\mbox{\scriptsize diff}} \to \infty$ when $T \stackrel{>}{\to} T_c$. Therefore, the time variation of $\rho$ is faster as $T-T_c$ gets smaller. We can rescale time to a faster time-scale by setting $t' = D_{\mbox{\scriptsize diff}} \, t = \frac{TT_c}{T-T_c} t$. Then, in the rescaled variables, the diffusion equation (\ref{eq:heat}) becomes independent of $T-T_c$ and gives the standard heat equation \,  $\partial_t \rho - \Delta_x \rho = 0$ \, (omitting the primes for simplicity).
\label{rem_heat_fast}
\end{remark}

\subsubsection{Case $T<T_c$ (small noise)}
\label{subsubsec:small_noise}

Let us now consider the dynamics of the system in the case of smaller noise.  In this section, we aim to prove the following:

\begin{proposition}
Let $T<T_c$ and define $c^2 = a^2-(d+2)T = (d+2) (T_c - T) >0$. Assume that the solution $(\rho^\tau, u^\tau)$ of system (\ref{E:alpha1_del}), (\ref{E:alpha2_del}) is smooth and converges smoothly towards a pair $(\rho,u)$ as $\tau \to 0$. Assume additionally that $(\rho^\tau,u^\tau)$ is not identically equal to $(\rho_0,0)$ where $\rho_0$ is constant in both space and time. 
Then, $u = c \omega$, where $\omega \in {\mathbb S}^{d-1}$ and the pair $(\rho,\omega)$ satisfies the following system:
\begin{align}
&\partial_t \rho + \nabla_x \cdot (c \rho \omega ) = 0 , \label{E:alpha1_st}\\
&(\partial_t + c \omega \cdot \nabla_x) \omega + \frac{T}{c} P (\nabla_x \ln \rho ) = 0 ,  \label{E:alpha2_st}
\end{align}
where $P = \mbox{Id} - \omega \otimes \omega$ is the orthogonal projection onto the hyperplane orthogonal to $\omega$. 
\label{prop_small_noise}
\end{proposition}

\begin{remark}
We easily see, taking the scalar product of (\ref{E:alpha2_st}) with $\omega$, that 
\begin{align*}
(\partial_t + c \omega \cdot \nabla_x) |\omega|^2  = 0 . 
\end{align*}
This implies that $|\omega(x,t)| \equiv 1$ for all time provided that $|\omega(x,0)| = 1$ initially. 
\label{rem_omega_norm1}
\end{remark}

\medskip
\noindent
{\bf Proof.} Equation \eqref{E:alpha2_del} is now written as:
\begin{equation} \label{E:Z}
(\partial_t + u^\tau \cdot \nabla_x ) u^\tau + T \nabla_x \ln \rho^\tau = \frac{1}{\tau a^2} u^\tau \left( c^2 - |u^\tau|^2 \right) .
\end{equation}
We now look for the equilibria of the right-hand side of (\ref{E:Z}), namely the solutions of 
$$  u \left( c^2 - |u|^2 \right) = 0. $$
There are two sets of equilibria. The first set reduces to the single point $u=0$. The second set is the sphere $|u|^2 = c^2$. We now show that in the spatially homogeneous setting, the first equilibrium is unstable, while the second class of equilibria is orbitally stable. 
Indeed, let us consider the solution $u(t)$ of the differential equation
$$ \frac{du}{dt} =  u (c^2 - |u|^2), \quad u(0) = u_0 .$$
Its solution can be analytically given by 
$$ |u(t)|^2 = \frac{|u_0|^2 c^2}{(c^2 - |u_0|^2) e^{-2 c^2 t} +  |u_0|^2}, $$
and is such that $|u(t)| \to c$ for all initial data $u_0$ except $u_0 = 0$. This shows that $u_0=0$ is an unstable equilibrium. On the other hand, if one perturbs an equilibrium $|u_0|^2 = c^2$ by a small amount, the solution will relax to an element of the circle $|u|^2 = c^2$ (where $u$ may be different from $u_0$).

Now, we let $\tau \to 0$ in the non spatially homogeneous system (\ref{E:alpha1_st}), (\ref{E:Z}). Unless $u^\tau$ is identically zero, which can only occur for a uniform density $\rho_0$, $u^\tau$ converges (at least formally) towards one of the stable equilibria. Therefore, we have $u^\tau \to u = c \omega$ with $|\omega| = 1$. In order to find the equation satisfied by $\omega$, we introduce a polar decomposition of the solution $u^\tau$. We write $u^\tau = c^\tau \omega^\tau$, with $c^\tau = |u^\tau|$ and $\omega^\tau = u^\tau/c^\tau$. Let $P^\tau$ the orthogonal projection of ${\mathbb R}^d$ onto the hyperplane orthogonal to $\omega^\tau$. The projection $P^\tau$ can be written tensorwise as $ P^\tau = \mbox{Id} - \omega^\tau \otimes \omega^\tau$. Inserting the polar decomposition of $u^\tau$ in (\ref{E:Z}) leads to:
\begin{equation} \label{E:Z2}
\omega^\tau (\partial_t + u^\tau \cdot \nabla_x ) c^\tau + c^\tau (\partial_t + u^\tau \cdot \nabla_x ) \omega^\tau + T \nabla_x \ln \rho^\tau = \frac{1}{\tau a^2} \omega^\tau \left( c^2 - (c^\tau)^2 \right) .
\end{equation}
We note that, because $|\omega^\tau|=1$ and the operator $\partial_t + u^\tau \cdot \nabla_x$ is a derivative, the vector $(\partial_t + u^\tau \cdot \nabla_x ) \omega^\tau$ is orthogonal to $\omega^\tau$. Furthermore, the first term of the left-hand side and the right-hand side are parallel to $\omega^\tau$. Therefore, applying $P^\tau$ to the second term of the left-hand side leaves it unchanged, while it cancels the first term of the left-hand side and the right-hand side. Consequently, applying $P^\tau$ to (\ref{E:Z2}) leads to:  
\begin{equation} \label{E:Z3}
c^\tau (\partial_t + u^\tau \cdot \nabla_x ) \omega^\tau + T P^\tau (\nabla_x \ln \rho^\tau ) = 0 .
\end{equation}
Now, taking the limit $\tau \to 0$ formally leads to (\ref{E:alpha2_st}). \endproof

\medskip
\noindent
System (\ref{E:alpha1_st}), (\ref{E:alpha2_st}) has been referred to in the literature as the Self-Organized Hydrodynamic (SOH) system \cite{Degond_etal_submitted1, Degond_etal_submitted2}. It has the form of a compressible gas dynamics system with isothermal equation-of-state and geometric constraint $|\omega|=1$. The projection operator $P$ which multiplies the pressure term $\nabla_x \ln \rho$ maintains the constraint over the course of time (see Remark \ref{rem_omega_norm1}). It results in a non-conservative term (since $P$ depends on $\omega$). This type of system has been derived for the first time in \cite{Degond_Motsch_M3AS08} and has been shown to be hyperbolic. Beyond this result, the mathematical study of such systems is in its infancy. A local existence result is given in \cite{Degond_etal_submitted2} and some special solutions are given in \cite{Motsch_Navoret_MMS11}. 

When $T \stackrel{<}{\to} T_c$, we have $c \to 0$ and therefore, $\rho$ becomes constant. To find a non-trivial dynamics, we must rescale time to a slower time-scale. By the time rescaling $t' = c t$, system (\ref{E:alpha1_st}), (\ref{E:alpha2_st}) can be written (dropping the primes): 
\begin{align*}
&\partial_t \rho + \nabla_x \cdot (\rho \omega ) = 0 , \\
&(\partial_t + \omega \cdot \nabla_x) \omega + \frac{T}{(d+2) (T_c - T)} P (\nabla_x \ln \rho ) = 0 ,  
\end{align*}
The parameter $(T_c-T)/T$ plays the role of the squared Mach-number in the standard compressible Euler system. Therefore, the limit $T \to T_c$, is similar to a small Mach-number limit. When $T \stackrel{<}{\to} T_c$, we formally get $\rho \to \rho_0$ where $\rho_0$ is a constant in space. With appropriate boundary conditions we can assume that $\rho_0$ is also independent of time. Then, the limit system as $T \stackrel{<}{\to} T_c$ is written as follows: 
\begin{align*}
&\nabla_x \cdot \omega = 0 , \\
&(\partial_t + \omega \cdot \nabla_x) \omega +  P \nabla_x \pi = 0 ,  
\end{align*}
where $\pi$ is a hydrodystatic pressure defined by 
$$\pi = \lim_{T \stackrel{<}{\to} T_c} \left( \frac{T}{(d+2) (T_c - T)} (\ln \rho - \ln \rho_0) \right) .$$
This system has been already proposed in \cite{Degond_etal_JSP10} for modeling gregariousness. It is a system of incompressible Euler equations subject to the geometric constraint $|\omega|=1$. 

The behavior of the (appropriately rescaled) SOH system (\ref{E:alpha1_st}), (\ref{E:alpha2_st}) when $T \stackrel{<}{\to} T_c$ is very different from the behavior of the (appropriately rescaled) diffusion equation (\ref{eq:heat}) when $T \stackrel{>}{\to} T_c$. The former is given by the incompressible Euler equations with geometric constraint $|\omega|=1$ while the latter is the standard diffusion equation (see Remark \ref{rem_heat_fast}). This drastic change of type is the signature of a phase transition further elaborated on in Section \ref{subsubsec:comments}.

\subsubsection{Comments}
\label{subsubsec:comments}

The fast relaxation limit of the compressible Euler system with self-propulsion (\ref{E:alpha1}), (\ref{E:alpha2}) exhibits two  different regimes for temperatures below or above the critical temperature $T_c$. Above $T_c$, the behavior of the system is that of a diffusion, while below $T_c$, the system obeys a hyperbolic system, the Self-Organized Hydrodynamic (SOH) model (\ref{E:alpha1_st}), (\ref{E:alpha2_st}). Even when the temperature is close to the critical temperature, it is not possible to match the two types of models in a smooth way, as noticed in the previous section. This abrupt change in the type of the model as the temperature crosses a threshold is a manifestation of a phase transition. Since $T_c$ depends on $a$, the phase transition originates from the self-propulsion force. Indeed, when $a=0$ (i.e. no self-propulsion), $T_c=0$ and there is no phase transition: the system is in the diffusion regime in all circumstances. 

Phase transitions in self-propelled particle systems have already been evidenced \cite{Vicsek_etal_PRL95} and an abundant literature has been devoted to them (see the review \cite{Vicsek_Zafeiris}). Phase transition in hydrodynamic models of self-propelled particles have first been studied in \cite{Toner_Tu_PRL95} (see also the review \cite{Toner_etal_AnnPhys05}). However, the equations proposed in \cite{Toner_Tu_PRL95} are more complicated than the ones seen here. They are derived solely on heuristic principles and invariance considerations. Their analysis is based on a combination of linear and nonlinear techniques. There is no link to the underlying particle models. The link between hydrodynamic and particle models of self-propeled particle systems has been made in \cite{Bertin_etal_PRE06, Bertin_etal_JPA09}, but for binary interaction mechanisms instead of the mean-field interaction considered here. Here, the hydrodynamic model is derived from the underlying particle dynamics and is much simpler: it merely consists of the isothermal compressible Euler model complemented with the self-propulsion force. The phase transition manifests itself in the change of type of the PDE which describes the system under large self-propulsion. 

\begin{figure}
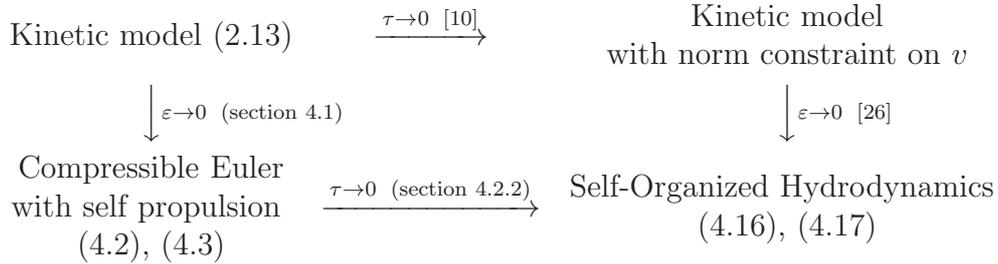

\begin{center}
$$ 
\begin{CD} 
\label{CD}
\mbox{Kinetic model (\ref{E:ourPDE_1})} 
@>{\tau \to 0  \, \mbox{  \scriptsize  \cite{Bostan_Carrillo}}  }>>  
\begin{array}{c} \mbox{Kinetic model} 
\\ \mbox{with norm constraint on $v$}
\end{array}
\\
@VV{\varepsilon \to 0 \, \mbox{ \scriptsize (section \ref{subsec_hydro})}  }V  @VV{\varepsilon \to 0 \, \mbox{  \scriptsize  \cite{Degond_Motsch_M3AS08}} } V\\
\begin{array}{c} 
\mbox{Compressible Euler} \\
\mbox{with self propulsion } \\
\mbox{(\ref{E:alpha1}), (\ref{E:alpha2})} 
\end{array}
@>{\tau \to 0 \, \mbox{ \scriptsize (section \ref{subsubsec:small_noise})}  }>>
\begin{array}{c} 
\mbox{Self-Organized Hydrodynamics } \\
\mbox{(\ref{E:alpha1_st}), (\ref{E:alpha2_st})} 
\end{array}
\end{CD}
$$
\caption{Two strategies. First strategy:  first take the relaxation limit $\tau \to 0$ (top horizontal arrow) and then pass to the hydrodynamic limit $\varepsilon \to 0$ in the resulting constrained kinetic model (right vertical arrow). This is done in \cite{Degond_Motsch_M3AS08} and \cite{Bostan_Carrillo}. Second strategy: first pass to the hydrodynamic limit $\varepsilon \to 0$ in the original kinetic model (left vertical arrow) and then pass to the fast relaxation limit $\tau \to 0$ in the resulting Euler system (bottom horizontal arrow). This is what is done here. }
\label{fig:CD}
\end{center}
\end{figure}

The SOH model (\ref{E:alpha1_st}), (\ref{E:alpha2_st}) has previously been derived in \cite{Degond_Motsch_M3AS08} from a system of self-propelled particles which have constant and uniform velocity. In \cite{Bostan_Carrillo}, it has been shown that the kinetic model with the velocity norm constraint of \cite{Degond_Motsch_M3AS08} is the the fast relaxation limit $\tau \to 0$ of (\ref{E:ourPDE_1}). A natural question is then whether the imposition of the norm constraint at the particle and kinetic levels as in \cite{Degond_Motsch_M3AS08} is necessary or useful. Indeed, there are now two ways of deriving the SOH model from (\ref{E:ourPDE_1}), which are summarized in Fig. \ref{fig:CD}. The first way is to follow the top horizontal and right vertical arrows successively. This is what is done in \cite{Degond_Motsch_M3AS08} and \cite{Bostan_Carrillo}. The second way is to follows the left vertical and bottom horizontal arrows successively. This is what is done here. 

There is however one noticeable difference between the two strategies. In \cite{Degond_Motsch_M3AS08}, the system is written: 
\begin{align}
&\partial_t \rho + \nabla_x \cdot (c_1 \rho \omega ) = 0 , \label{eq:SOH_rho}\\
&(\partial_t + c_2 \omega \cdot \nabla_x) \omega + \delta P (\nabla_x \ln \rho ) = 0 , \label{eq:SOH_omega}
\end{align}
where the constants $c_1$ and $c_2$ are such that $0 < c_2 \leq c_1$. In \cite{Frouvelle_M3AS12}, it has even been shown that taking into account an anisotropic vision, ($c_1, c_2$) can be made arbitrary. We would recover (\ref{E:alpha1_st}), (\ref{E:alpha2_st}) if we could set $c_2 = c_1 = c$ and $\delta = T/c$. However, it is not possible since, in \cite{Degond_Motsch_M3AS08}, $c_2 \not = c_1$. Therefore, the diagram in Fig. \ref{fig:CD} is {\bf not} commutative. This shows how important it is to enforce the norm constraint at the kinetic level, in spite of the induced difficulties. Indeed, such constrained kinetic models do not exhibit momentum conservation and lack classical collision invariants such as those derived in Section \ref{subsec:CI}. This makes the derivation of hydrodynamic models considerably more complex. In \cite{Degond_Motsch_M3AS08}, new concepts have been developed to bypass these difficulties. The fact that $c_1 \not = c_2$ is a direct consequence of these features. 

The question whether $c_1 = c_2$ or not has important consequences. Indeed, the SOH system is Galilean invariant if and only if $c_1 = c_2$ (which is the case here). However, in the generic case $c_1 \not = c_2$,  the SOH system is {\bf not} Galilean invariant.  This fact reflects a key feature of collective motion: the anisotropy of information flow. This is evidenced in car traffic, where perturbations of velocities (typically moving with speed $c_2<c_1$) propagate upstream the flow (whose speed is $c_1$). This phenomenon is a consequence of the fact that information propagates upstream from drivers ahead to drivers behind.  When $c_1 = c_2$, this property is lost and information spreads in an isotropic way just like usual gas dynamics. In such a case, the model is unable to correctly reproduce the complex emerging patterns (such as congestions and waves in car traffic). We refer to \cite{Tu_etal_PRL98} for a discussion of the loss of Galilean invariance in biological swarms. 

The approach of \cite{Degond_Motsch_M3AS08} does not lead to phase transitions: whatever the values of the temperature $T$ or self-propulsion velocity $a$, the limit system is of hydrodynamic type. The very large intensity of the self propulsion force which is needed to perform the $\tau \to 0$ limit first (top horizontal arrow in diagram \ref{fig:CD}) prevents any diffusion regime to establish. However, the emergence of a phase transition is a well-established experimental fact (see review in \cite{Vicsek_Zafeiris}). The inability of \cite{Degond_Motsch_M3AS08} to produce phase transitions can be seen as a major drawback. Fortunately, it has been proved in \cite{Degond_etal_submitted1} that adding a dependence of the social force intensity $\sigma^{-1}$ upon the flux $j$ restores the phase transition. It also brings density-dependent phase transition. More precisely, the parameter controling the phase transition is the ratio $T/\rho$. A natural question is to explore similar features here. This question will be investigated in future work.

Here, we choose a different direction: we explore whether the inclusion of diffusion terms can cure the deficiency of the SOH models derived here (i.e. the fact that $c_1=c_2$). This is the goal of the next section.

\section{Diffusive corrections to the hydrodynamic model}
\label{sec:diffusive}

\subsection{Setting of the problem}
\label{subsec:expan}

In this section, we derive ${\mathcal O}(\varepsilon)$ diffusive corrections to the compressible Euler model with self propulsion (\ref{E:alpha1}), (\ref{E:alpha2}). We show that these diffusive corrections lead to a system of isothermal compressible Navier-Stokes equations including self-propulsion. In a second step, we perform the fast relaxation limit $\tau \to 0$ in the resulting Navier-Stokes model and compare it to the SOH model obtained in Proposition \ref{prop_small_noise}. 

We now set up the problem.  We define 
\begin{equation} 
\rho_f = \int f \, dv, \quad j_f = \int f \, v \, dv = \rho_f u_f. 
\label{eq:momf}
\end{equation}
We assume that the observation kernel $K(|\xi|)$ is such that $\int K(|\xi|) \, d \xi = 1$ and we denote by $k >0$ the second moment of $K$, i.e.
$$ \frac{1}{2} \int K(|\xi|) \, (\xi \otimes \xi) \, d \xi = k \, \mbox{Id}. $$ 
For instance, if $K$ is the indicator function of the ball of radius $1$, $k = \frac{|{\mathbb S}^{d-1}|}{2 d (d+2)}$, where $|{\mathbb S}^{d-1}|$ is the $(d-1)$-dimensional measure of ${\mathbb S}^{d-1}$. In the cases of $d=2$ and $d=3$, we respectively get $k = \pi/8$ and $k = 2 \pi / 15$. 

We first give the expansion of ${\mathcal F}_f^\varepsilon$ up to the fourth order in $\varepsilon$.

\begin{lemma}
We have:
\begin{eqnarray}
&& \bar v_f^\varepsilon = u_f + \varepsilon^2 u_f^1 + O(\varepsilon^4), \label{eq:vfeps2}\\
&& {\mathcal F}_f^\varepsilon = \sigma^{-1} (u_f - v) + \varepsilon^2 \sigma^{-1} u_f^1 + O(\varepsilon^4), \label{eq:Ffeps2}
\end{eqnarray}
where
\begin{eqnarray}
&& u_f^1 = \frac{k_R}{\rho_f^2} \, \big(\rho_f \Delta_x j_f - j_f \Delta_x \rho_f \big) = k_R \, \big(\Delta_x u_f + 2 (\nabla_x \ln  \rho_f \cdot \nabla_x) u_f \big) , \label{eq:uf1}
\end{eqnarray}
and $k_R = k R^2$. 
\label{lem:expanFf}
\end{lemma}

\medskip
\noindent
We give the proof of this very simple lemma in Appendix \ref{app_expan}. The $O(\varepsilon^2)$ correction to the mean velocity $u_f$ takes into account the non-local character of the average (\ref{vfeps}). This correction involves gradients of the local density and flux. They quantify how information spreads due to the fact that the agents observe their environment over a certain spatial extent. Because the observation is supposed isotropic, the correction only involves $O(\varepsilon^2)$ terms and second order derivatives. In the case of non-isotropic observation kernels, $O(\varepsilon)$ corrections involving first order gradients would be obtained. The study of this effect is postponed to future work. 

Now, up to terms of order $\varepsilon^4$ (which are dropped), equation (\ref{eq:rescaledkin}) is written: 
\begin{eqnarray} 
&&\hspace{-1cm}\veps \left[ \partial_t f^\veps + \nabla_x \cdot (v f^\veps) + \tau^{-1} \nabla_v \cdot \Big( \big(1-\frac{|v|^2}{a^2} \big) v f^\veps \Big) \right] \nonumber \\
&&\hspace{5cm} + \, \varepsilon^2 \sigma^{-1} \nabla_v \cdot (u_{f^\varepsilon}^1 f^\varepsilon) = Q(f^\veps), \label{E:ourPDE_10}
\end{eqnarray}
where again, $Q$ is given by (\ref{E:collisionAlphaZero}). Now, we look for a fluid model approximating (\ref{E:ourPDE_10}) which includes the $O(\varepsilon)$ terms. We will adopt Chapman-Enskog's method consisting in writing the $O(\varepsilon)$ terms in terms of spatial derivatives only (see e.g. the review in \cite{Degond_Birkhauser04}). The model is introduced and discussed in the next section.

\subsection{Compressible Navier-Stokes equations with self-propulsion}
\label{subsec_NS_statement}

The compressible Navier-Stokes system with self-propulsion is established in the following: 

\begin{theorem}
Let $f^\varepsilon$ be the solution of (\ref{E:ourPDE_10}) associated to a given initial condition $f_I$, which is supposed independent of $\varepsilon$ for simplicity. Let $(\rho_{f^\varepsilon}, \rho_{f^\varepsilon} u_{f^\varepsilon})$ be the moments of $f^\varepsilon$ defined by (\ref{eq:momf}). Then, we can formally write $(\rho_{f^\varepsilon}, \rho_{f^\varepsilon} u_{f^\varepsilon}) = (\rho^\varepsilon, \rho^\varepsilon u^\varepsilon) + O(\varepsilon^2)$, provided that $(\rho^\varepsilon, \rho^\varepsilon u^\varepsilon)$ satisfy the following set of compressible Navier-Stokes-like equations (where we drop the superscripts $\varepsilon$ upon $\rho^\varepsilon$ and $u^\varepsilon$ for the sake of clarity): 
\begin{align}
&\partial_t \rho + \nabla_x \cdot (\rho u ) = 0 , \label{E:alpha1_NS}\\
&\partial_t (\rho u ) + \lambda^\varepsilon \nabla_x \cdot (\rho u \otimes u ) + \nabla_x \pi^\varepsilon(\rho,|u|) = -\frac{1}{\tau^\varepsilon} \rho u \big(\frac{|u|^2}{a^2} - \chi^\varepsilon \big) \nonumber \\
& \hspace{2cm} + \varepsilon \big( \mu \nabla_x \cdot (\rho {\mathcal E}(u)) + \frac{k_R}{\sigma} \rho \Delta_x u \big) + \frac{2 \varepsilon k_R}{\sigma}(\nabla_x \rho \cdot \nabla_x) u  \nonumber \\
& \hspace{2cm} + \frac{\varepsilon \lambda}{2} \rho \big[ (\nabla_x \cdot u) u + \nabla_x \big( \frac{|u|^2}{2} \big) + (u \cdot \nabla_x) u \big] , \label{E:alpha2_NS}
\end{align}
where 
\begin{eqnarray}
&& \lambda^\varepsilon = 1 - \varepsilon \lambda, \quad \lambda = \frac{2 \sigma T}{\tau a^2}, \label{eq:lambdaeps} \\
&& \pi^\varepsilon(\rho,|u|) = T \rho - \varepsilon \pi(\rho,|u|), \quad \pi(\rho, |u|) = \frac{\lambda}{2} \rho \big\{ (d+2) T  - a^2  +  |u|^2 \big\}, \label{eq:pieps} \\
&& \chi^\varepsilon = \frac{1 -  \frac{\varepsilon \lambda}{2} (d+2) - (d+2)  \frac{T}{a^2} \big( 1 - \frac{\varepsilon \lambda}{2} (d+4) \big) }{1 - \frac{\varepsilon \lambda}{2} (d+8)}, \label{eq:chieps} \\
&&  \frac{1}{\tau^\varepsilon} = \frac{1}{\tau} \big( 1 - \varepsilon \frac{\lambda}{2} (d+8) \big), \label{eq:taueps} \\
& & \mu = \sigma T, \quad {\mathcal E}(u) = \frac{1}{2} (\nabla_x u + (\nabla_x u)^T). \label{eq:muE(u)}
\end{eqnarray}
\label{thm:NS}
\end{theorem}

\noindent
The proof of this result is fairly technical. It is given in Appendix \ref{app_NS_proof}. The coefficient $\mu$ is the fluid viscosity, while $\pi^\varepsilon(\rho,u)$ is a velocity-dependent pressure. 

\medskip
We compare this Navier-Stokes model to the hydrodynamic model (\ref{E:alpha1}), (\ref{E:alpha2}) and provide a physical interpretation of the ${\mathcal O}(\varepsilon)$ correction terms. The mass conservation equation (\ref{E:alpha1_NS}) is unchanged  compared to (\ref{E:alpha1}), as it should be. Indeed, the particle density is still a conservative variable transported by the fluid velocity. By contrast, there is a wealth of new terms in the momentum balance equation (\ref{E:alpha2_NS}). We can identify the only term which come from classical fluid viscosity: this is $\varepsilon \mu \nabla_x \cdot (\rho {\mathcal E}(u))$. The present viscous term is different from that appearing in the full Navier-Stokes system, where the temperature is determined by the energy balance equation. Here, the viscous term involves the symmetrized velocity gradient tensor ${\mathcal E}(u)$ instead of the rate of strain tensor ${\mathcal S}(u) = {\mathcal E}(u) - \frac{1}{d} (\nabla_x \cdot u) \mbox{Id}$. This difference is solely due to the isothermal character of the model and not to self-propulsion.

Then, there are two terms which come from the non-locality of the interaction and which are proportional to $k_R$. The first one contributes to adding more viscosity and is equal to $\varepsilon \frac{k_R}{\sigma} \rho \Delta_x u$. The second one contributes to convecting the velocity in the direction of the gradient of $\rho$. It is written $\frac{2 \varepsilon k_R}{\sigma}(\nabla_x \rho \cdot \nabla_x) u $. 

All other ${\mathcal O}(\varepsilon)$ terms originate from self-propulsion. The self-propulsion first contributes to a similar relaxation source term as in (\ref{E:alpha2}): the term $-\frac{1}{\tau^\varepsilon} \rho u \big(\frac{|u|^2}{a^2} - \chi^\varepsilon \big)$. However, both the relaxation rate $\frac{1}{\tau^\varepsilon}$ and the bulk comfort velocity $\sqrt{\chi^\varepsilon} a$ towards which this source term is relaxing are different from those of (\ref{E:alpha2}), respectively equal to $\frac{1}{\tau}$ and $\sqrt{a^2 - (d+2)T}$. We notice that the former are ${\mathcal O}(\varepsilon)$ corrections of the latter, as they should be. However, there is a second source term, the term 
$$\frac{\varepsilon \lambda}{2} \rho \big[ (\nabla_x \cdot u) u + \nabla_x \big( \frac{|u|^2}{2} \big) + (u \cdot \nabla_x) u \big],$$ 
which describes a force comprising three terms. The first term is just a friction proportional to the compressibility $\nabla_x \cdot u$. The second term is a force acting in the direction of the gradient of $|u|^2$. Finally, the last term is a force acting in the direction of $(u \cdot \nabla_x) u$. 

The self-propulsion force also induces some changes in the convection terms. The first one is a change in the convection velocity of the momentum, which is not just $u$ but $\lambda^\varepsilon u$. The coefficient $\lambda^\varepsilon$ is close to one but not exactly equal to one, the difference being ${\mathcal O}(\varepsilon)$. This feature makes the model closer to the generic SOH model (\ref{eq:SOH_rho}), (\ref{eq:SOH_omega}), which has different convection velocities $c_1$ and $c_2$ for the density and velocity. The last contribution of the self-propulsion is a modification of the pressure. The isothermal pressure $\rho T$ is complemented by an ${\mathcal O}(\varepsilon)$ pressure correction $\varepsilon \pi(\rho,|u|)$ which depends on the norm of the velocity. This correction is positive or negative according to whether the speed is larger or smaller than the bulk comfort velocity $\sqrt{a^2 - (d+2)T}$ of the non-viscous model. 

We now discuss this model in view of the model proposed in \cite{Toner_Tu_PRL95, Toner_etal_AnnPhys05}, which is considered as the `paradigmatic' 
model of hydrodynamic type for flocking. In this model, the mass conservation equation is the same as (\ref{E:alpha1_NS}) but the momentum balance equation takes the following form (in the absence of extrernal force): 
\begin{eqnarray}
&& \hspace{-1cm} 
\partial_t u  + (u \cdot \nabla_x) u = \alpha u - \beta |u|^2u - \nabla_x P + D_L \nabla_x (\nabla_x \cdot u) \nonumber \\
&& \hspace{5cm} \, +  D_1 \Delta_x u  + D_2 (u \nabla_x \cdot u)^2 u , 
\label{eq:vel_TT}
\end{eqnarray}
where $\alpha$, $\beta$, $D_L$, $D_1$ and $D_2$ are positive coefficients and $P =P(\rho)$ is the pressure, which depends nonlinearly of $\rho$. This model has been derived on the basis of invariance considerations. In order to compare our model with (\ref{eq:vel_TT}), we use the mass conservation equation (\ref{E:alpha1_NS}) to write the momentum balance equation (\ref{E:alpha2_NS}) as follows: 
\begin{eqnarray}
&& \hspace{-1cm} 
\partial_t u  + \big( 1 - \varepsilon \frac{3 \lambda}{2} \big) \, (u \cdot \nabla_x) u - \varepsilon \lambda \, (u \cdot \nabla_x \ln \rho) \, u  - \varepsilon \frac{3 \lambda}{2} \, (\nabla_x \cdot u)  \, u - \varepsilon \frac{\lambda}{2} \, \nabla_x \big( \frac{|u|^2}{2} \big) \nonumber \\
&& \hspace{0.5cm} = -\frac{1}{\tau^\varepsilon} \rho u \big(\frac{|u|^2}{a^2} - \chi^\varepsilon \big) - \nabla_x \pi^\varepsilon(\rho,|u|)
+ \varepsilon \big( \mu + \frac{2 k_R}{\sigma} \big) (\nabla_x \ln \rho \cdot \nabla_x) u 
\nonumber \\
&& \hspace{2cm} 
+ \varepsilon \big( \mu + \frac{k_R}{\sigma} \big) \Delta_x u + \varepsilon \mu (\nabla_x u) \nabla_x \ln \rho + \varepsilon \mu \nabla_x (\nabla_x \cdot u). 
\label{eq:vel_BD}
\end{eqnarray}
There is one term in (\ref{eq:vel_TT}) which is missing from (\ref{eq:vel_BD}): the anisotropic velocity diffusion $D_2 (u \cdot \nabla_x)^2 u $. By contrast, there are many terms appearing in (\ref{eq:vel_BD}) which are not present in (\ref{eq:vel_TT}): the third, fourth and fifth terms of the left-hand side and the third and fifth terms of the right-hand side. Additionally, among the terms which are common to both formulas (namely the first and second terms of the left-hand side and the first, second, fourth and sixth terms of the right-hand side of (\ref{eq:vel_BD})), some of them assume different forms. Indeed, the second term (the convection term $(u \cdot \nabla_x) u$) is multiplied by the constant $\big( 1 - \varepsilon \frac{3 \lambda}{2} \big)$ less than one in (\ref{eq:vel_BD}), while this coefficient is exactly one in (\ref{eq:vel_TT}). This difference is significant, in view of the previous discussion about mass and momentum convection velocities in swarming systems. Another difference is in the pressure term (second term at the right-hand side of (\ref{eq:vel_BD})). In our model, the pressure depends on both the density and the norm of the velocity, while it depens on the density only in (\ref{eq:vel_TT}). By contrast, the dependence upon the density is linear in our case, while it is nonlinear in (\ref{eq:vel_TT}). This discussion illustrates that phenomenological models can differ significantly from first principle models when complex phenomena such as swarming behavior are concerned. The question whether these differences lead to perceivable changes in the qualitative behavior of the solution has not yet been investigated.

The mathematical properties of this system (such as e.g. the stability of this system for small perturbations of a homogeneous state) will be studied in future work. In this work, we investigate the fast relaxation limit $\tau \to 0$. This is the goal of the next section.

\subsection{The fast relaxation limit in the compressible Navier-Stokes equations with self-propulsion}
\label{subsec_NS_fast}

In this section, we examine the limit $\tau \to 0$ in the Navier-Stokes system (\ref{E:alpha1_NS}), (\ref{E:alpha2_NS}). Since this system was derived under the assumption that $\varepsilon$ is small, $\varepsilon$ needs to tend to $0$ at least as fast as $\tau$ tends to $0$. Here, we decide to make $\varepsilon$ and $\tau$ proportional. This is the borderline case, because, to be consistent, we should have linked $\varepsilon$ and $\tau$ in this way already at the kinetic level. However, as pointed out earlier, the analysis of this scaling limit is more complex and is still under scrutiny \cite{Barbaro_etal_inprep}. Our conjecture is that the kind of model we get is the same in both limits. Therefore, investigating it at the level of the Navier-Stokes system is a preparation before performing the limit directly from the kinetic level.

Since $\lambda$ is proportional to $1/\tau$, we decide to relate $\varepsilon$ and $\lambda$ in such a way that $\varepsilon \lambda$ is a constant $\alpha$, i.e.
\begin{eqnarray}
&& \alpha = \varepsilon \lambda. \label{eq:al}
\end{eqnarray}
We need to keep in mind that, strictly speaking, $\alpha$ must be $\ll 1$ otherwise the derivation of the Navier-Stokes model in the previous section loses its valididty. We decide to express $\varepsilon$ as a function of $\tau$. We can write
\begin{eqnarray}
&& \varepsilon = \kappa_\alpha \tau, \quad \quad \kappa_\alpha = \frac{\alpha a^2}{2 \sigma T}. \label{eq:kappa}
\end{eqnarray}
We also define 
\begin{eqnarray}
&& \lambda_\alpha := \lambda^{\frac{\alpha}{\lambda}} = 1 - \alpha, \label{eq:lambda_al} \\
&& \pi_\alpha(\rho,|u|) := \pi^{\frac{\alpha}{\lambda}}(\rho,|u|) =T \rho - \frac{\alpha}{2} \rho \big\{ (d+2) T  - a^2  +  |u|^2 \big\}, \label{eq:pi_al} \\
&&  \frac{1}{\tau_\alpha} = \frac{1}{\tau^{\frac{\alpha}{\lambda}}} = \frac{\xi_\alpha}{\tau}, \quad \quad \xi_\alpha =   1 -  \frac{d+8}{2} \alpha . \label{eq:tau_al}
\end{eqnarray}
We can introduce an $\alpha$-dependent critical temperature $T_c(\alpha)$ by: 
\begin{eqnarray}
&& T_c(\alpha) = T_c(0) \,  \frac{1 - \frac{d+2}{2}\alpha}{1 - \frac{d+4}{2}\alpha},  \label{eq:Tc(alpha)}
\end{eqnarray}
where $T_c(0)$ is the critical temperature (\ref{eq:Tth}). Then, we can write the square of the comfort velocity in the Navier-Stokes model $c_1(\alpha)$ as 
\begin{eqnarray}
&& c_1^2 (\alpha) = a^2 \chi^{\frac{\alpha}{\lambda}} = (d+2) \, \frac{1 - \frac{d+2}{2}\alpha}{1 - \frac{d+4}{2}\alpha} \,  (T_c(\alpha) - T).  \label{eq:ch_al} 
\end{eqnarray}
Since $\alpha$ must be small, we can limit its range in such a way that $\xi_\alpha >0$ and that the relaxation time $\tau_\alpha$ remains positive. Therefore, we have $\alpha \in [0,\frac{2}{d+8}]$. In this interval, $T_c(\alpha)$ is an increasing function of $\alpha$ with values in $[T_c(0), \frac{3}{2}T_c(0)]$. We also assume that the temperature is below the critical temperature $T<T_c(\alpha)$ in such a way that $c_1^2 (\alpha) >0$. 

We finally denote by $(\rho^\tau, u^\tau)$ the solution of the Navier-Stokes system (\ref{E:alpha1_NS}), (\ref{E:alpha2_NS}). We re-write the system in the new notation: 
\begin{align}
&\partial_t \rho^\tau + \nabla_x \cdot (\rho^\tau u^\tau ) = 0 , \label{E:alpha1_NS2}\\
&\partial_t (\rho^\tau u^\tau ) + \lambda_\alpha \nabla_x \cdot (\rho^\tau u^\tau \otimes u^\tau ) + \nabla_x \pi_\alpha(\rho^\tau,|u^\tau|) = -\frac{1}{\tau} \, \frac{\xi_\alpha}{a^2} \rho^\tau u^\tau \big(|u|^2 - c_1^2(\alpha) \big) \nonumber \\
& \hspace{2cm} + \tau \, \frac{1}{\kappa_\alpha} \left\{  \big( \mu \nabla_x \cdot (\rho^\tau {\mathcal E}(u^\tau)) + \frac{k_R}{\sigma} \rho^\tau \Delta_x u^\tau \big) + \frac{2 k_R}{\sigma}(\nabla_x \rho^\tau \cdot \nabla_x) u^\tau \right\}  \nonumber \\
& \hspace{2cm} + \frac{\alpha}{2} \rho \big[ (\nabla_x \cdot u^\tau) u^\tau + \nabla_x \big( \frac{|u^\tau|^2}{2} \big) + (u^\tau \cdot \nabla_x) u^\tau \big] , \label{E:alpha2_NS2}
\end{align}
Now, we formally let $\tau \to 0$ in this system, keeping all other parameters fixed, and in particular $\alpha$. We get the following:

\begin{proposition}
Let $\alpha \in [0,\frac{2}{d+8}]$ and let $T<T_c(\alpha)$. Assume that the solution $(\rho^\tau, u^\tau)$ of system (\ref{E:alpha1_NS2}), (\ref{E:alpha2_NS2}) is smooth and converges smoothly towards a pair $(\rho,u)$ as $\tau \to 0$. Assume additionally that $(\rho^\tau,u^\tau)$ is not identically equal to $(\rho_0,0)$ where $\rho_0$ is constant in both space and time. 
Then, $u = c_1(\alpha) \omega$, where $\omega \in {\mathbb S}^{d-1}$ and the pair $(\rho,\omega)$ satisfies the following system:
\begin{align}
&\partial_t \rho + \nabla_x \cdot (c_1(\alpha) \rho \omega ) = 0 , \label{E:alpha1_stdif}\\
&(\partial_t + c_2(\alpha) \omega \cdot \nabla_x) \omega + \delta_\alpha P (\nabla_x \ln \rho ) = 0 ,  \label{E:alpha2_stdif}
\end{align}
where $P = \mbox{Id} - \omega \otimes \omega$ is the orthogonal projection onto the hyperplane orthogonal to $\omega$ and 
\begin{eqnarray}
&& c_2 (\alpha) = \big( 1 - \frac{3}{2} \alpha \big) c_1 (\alpha), \label{eq:c2_al} \\
&& \delta_\alpha = \frac{T_\alpha}{c_1 (\alpha)}, \quad \quad T_\alpha = \frac{1}{1 - \frac{d+8}{2}\alpha} \left\{ \left( 1 +  \frac{d-4}{2} \alpha\right) T - \frac{3}{2} \alpha a^2 \right\}.  \label{eq:del_al} 
\end{eqnarray}
\label{prop_small_noise_diff}
\end{proposition}

\begin{remark}
For small $\alpha$, up to terms of order $\alpha^2$, we have:
\begin{eqnarray*} 
&& c_1(\alpha) = c_1(0) + \alpha (d+2) (3 T_c(0) - 2T ) > c_1(0), \\
&& c_2(\alpha ) = c_1(0) + \alpha \frac{d+2}{2} (3 T_c(0) - T ) > c_2(0). 
\end{eqnarray*} 
Therefore, the convection speeds are larger when first order corrections are included. Additionally, it is easy to prove that $c_1(\alpha)$ is an increasing function of $\alpha$ on its interval of definition.
\label{rem_alpha_to_0}
\end{remark}

\medskip
\noindent
{\bf Proof.} The proof is similar to that of Proposition \ref{prop_small_noise_diff} and is only sketched. First, observe that the second line of the right-hand side of (\ref{E:alpha2_NS2}), which is proportional to $\tau$, simply vanishes in the limit. Since all the diffusion terms in the Navier-Stokes system are contained on this line, there is no diffusion in the limit system. Second, since $|u^\tau| \to c_1(\alpha)$, which is a constant, the second term of the third line of (\ref{E:alpha2_NS2}) also vanishes in the limit. For the same reason, the pressure $\pi_\alpha(\rho^\tau,|u^\tau|) \to \pi_\alpha(\rho,c_1(\alpha)) = T_\alpha \rho$. Therefore, we recover an isothermal pressure equation-of-state in the limit. Then, since (\ref{eq:del_al}) is obtained by projecting (\ref{E:alpha2_NS2}) onto the hyperplane normal to $u^\tau$, the first term of the third line of (\ref{E:alpha2_NS2}), which is parallel to $u^\tau$, also vanishes in the limit. Finally, the last term of the third line of (\ref{E:alpha2_NS2}) combines with the second term at the left-hand side of (\ref{E:alpha2_NS2}) and yields the second term of the left-hand side of (\ref{E:alpha2_stdif}). Indeed, we get 
$$ P^\tau \big[ \lambda_\alpha \nabla_x \cdot (\rho^\tau u^\tau \otimes u^\tau ) - \frac{\alpha}{2} \rho (u^\tau \cdot \nabla_x) u^\tau \big] \to (\lambda_\alpha -  \frac{\alpha}{2}) \rho c_1^2(\alpha) (\omega \cdot \nabla_x) \omega, $$
where $P^\tau$ is defined like in the proof of Proposition \ref{prop_small_noise_diff}. Therefore, $c_2(\alpha) = (\lambda_\alpha -  \frac{\alpha}{2}) c_1(\alpha)$ is given by (\ref{eq:c2_al}). The remaining details are left to the reader. \endproof

Now, system (\ref{E:alpha1_stdif}), (\ref{E:alpha2_stdif}) is in the form of the generic SOH model (\ref{eq:SOH_rho}), (\ref{eq:SOH_omega}) with $c_2 < c_1$, like in \cite{Degond_Motsch_M3AS08}. We conclude that the inclusion of diffusive terms in the hydrodynamic model before passing to the fast relaxation limit was a successful approach. It has cured the deficiency of the SOH model derived in Section \ref{subsubsec:small_noise}. In some sense, the purely hydrodynamic model obtained at Theorem \ref{thm:hydro} is too simple to describe the full complexity of the system but the inclusion of diffusive terms is enough to restore the adequate level of complexity. 

It is intriguing that the only diffusive corrections that are kept in the fast relaxation limit are those coming from the self-propulsion force. It would be interesting to compute the ${\mathcal O}(\tau)$ corrections to this model. Then, the other diffusive terms (those arising from viscosity and non-locality of the interaction) would appear. The resulting model would then have more relevance as an approximation of the original kinetic model. It should also be compared to the diffusive SOH model obtained in \cite{Degond_Yang_M3AS10}, which is quite complex. It would be instructive to see if the present approach could help get a cleaner model.

\section{Conclusion and perspectives}
\label{sec_conclu}

In this work, we have provided evidence of a phase transition from disordered to ordered motion in a hydrodynamic model of socially interacting agents with self-propulsion. The model we have investigated has been derived from a particle system combining a noisy Cucker-Smale consensus force and self-propulsion. We have shown that the phase transition appears in the limit of a large self-propulsion force and manifests itself as a change of type of the limit model (from hyperbolic to diffusive) at the crossing of a critical noise intensity. We have also shown that, in the hyperbolic regime, the resulting SOH (self-Organized Hydrodynamics) model suffers from unnecessary restrictions on the range of its coefficients. To remove these restrictions, we have computed diffusive correction to the model. With these diffusive corrections, the restrictions on the SOH model obtained in the limit of a large propulsion force disappear. 

As pointed out in the core of the work, many points deserve further elucidation. A first one, currently under scrutiny, consists in performing the combined hydrodynamic and large self-propulsion force simultaneously at the level of the kinetic model. We anticipate that similar phase transitions will emerge and will be described by the same limit models, with possibly different coefficients. Other points would be worth being developed. The computation of diffusive corrections to the SOH model for instance would be of great practical use. Further investigations of quantities attached to the social force are also very promising, such as the role of a possible anisotropy of the observation kernels, or of a different scaling of its range. Finally, the understanding of the transition between the two phases and how the two different models can be matched at this transition is also crucial for applications.

\appendix
\section{Appendix: Properties of the collision operator $Q$}
\label{app_prop_Q}

{\bf Proof of Lemma \ref{lem:Q}:} Observe:
\begin{align*}
\nabla_v \cdot \left[ M_{u_f} \nabla_v \big( \frac{f}{M_{u_f}} \big) \right] &= \Delta_v f + \nabla_v \cdot \left[ -\nabla_v \left( \ln M_{u_f} \right) f \right]\\
&= \Delta_v f + \nabla_v \cdot \left[ \frac{v-u}{T} \, f\right].
\end{align*}
Then (\ref{E:collision2}) follows. Formula (\ref{eq:self}) is a consequence of Green's formula. Now, suppose that $f \in \mbox{Ker}\; Q$, then: 
\begin{align}
0 = \int Q(f) \frac{f}{M_{u_f}} \,dv = D \int M_{u_f} \Big|\nabla_v \big( \frac{f}{M_{u_f}} \big) \Big|^2 \,dv
\label{E:defOfF}
\end{align}
This implies that there exists a constant $C$ such that $f = C M_{u_f}$. In particular, $f$ is of the form $C M_u$ for a given vector $u \in {\mathbb R}^d$. Reciprocally, let $u \in \mathbb{R}^d$ and $C$ a constant, and construct $f = CM_u$.  Then 
\begin{align*}
u_f = \frac{\int f v \,dv}{\int f \,dv}
= \frac{\int M_u v \,dv}{ \int M_u \,dv}
= u. 
\end{align*}
Then $f = C M_{u_f}$. So, by \eqref{E:collision2}, $Q(f) = 0$ and thus $f \in \text{Ker}(Q)$.  This ends the proof of Lemma \ref{lem:Q}. \endproof

\medskip
\noindent
{\bf Proof of Proposition \ref{prop:CI}.} We first show that there are $d+1$ obvious linearly independent CI's. More precisely, suppose $\psi(v) =1$ or $\psi(v) =v$.  Then $\psi$ is a collision invariant. Indeed, clearly, such $\psi$ satisfies that $\psi \in V_{u_f}$ for all $f$ such that $f/M_{u_f} \in V_{u_f}$. The statement that $\psi(v) = 1$ is a CI follows from applying (\ref{eq:self0}) with $g=1$. Let now $k \in \{1, \ldots, d \}$ and consider  $\psi (v) = v_k$.  Then
\begin{align*}
\int Q(f) \, v_k \,dv &= - D \int M_{u_f} \partial_{v_k} \big(\frac{f}{M_{u_f}}\big) \,dv\\
&= - D \int \big( \partial_{v_k} f - f \partial_{v_k}(\ln M_{u_f}) \big) \,dv.
\end{align*}
Using Green's Theorem, the first term disappears, and it follows, from the definition of $u_f$,  that
\begin{align*}
\int Q(f) \, v_k \,dv &= -\frac{1}{\sigma} \int f\left(v_k - (u_{f})_k \right) \,dv= 0, 
\end{align*}
which shows that $\psi (v) = v_k$ is also a CI. Thus, ${\mathcal C}$ contains a $(d+1)$-dimensional linear space, namely $\text{Span}\{ 1, v_1, \dots, v_d \}$. We now prove that ${\mathcal C}$ is identically equal to this space. 

Let $u \in {\mathbb R}^d$. Then, we define the operator $R(u;f)$ as follows: 
\begin{equation} \label{E:collision_u}
R(u;f) = D \nabla_v \cdot \left[ M_{u} \nabla_v \big( \frac{f}{M_{u}} \big) \right],
\end{equation}
for all $f$ such that $f/M_{u} \in V_{u}$.
We notice that for given $u \in {\mathbb R}^d$, $R(u;f)$ is a linear operator with respect to $f$ and that 
$$ Q(f) = R(u_f;f).$$
Then, $\psi$ is a CI of $Q$ if and only if 
$$ \int  R(u_f;f) \psi \,dv = 0, \quad \forall f \quad \mbox{ s.t.} \quad f/M_{u_f} \in V_{u_f} \quad \mbox{ and } \quad \psi \in V_{u_f}, $$
or equivalently, if and only if:
\begin{eqnarray*} 
& & \hspace{-1cm} \forall u \in {\mathbb R}^d \quad \mbox{ we have: }  \quad \psi \in V_{u} \quad \mbox{ and } \quad \int  R(u;f) \psi \,dv = 0,  \\
& & \hspace{3cm}
\forall f \quad \mbox{such that} \quad f/M_{u} \in V_{u} \quad \mbox{ and } \quad  u_f = u .  
\end{eqnarray*}

As a first step, we fix $u \in {\mathbb R}^d$ and find all $\psi \in V_{u}$ which satisfy:
\begin{equation} \int  R(u;f) \psi \,dv = 0, \quad \forall f \quad \mbox{such that} \quad f/M_{u} \in V_{u} \quad \mbox{ and } \quad u_f = u . 
\label{eq:condpsi}
\end{equation}
Then, we will make $u$ arbitrary. We note that any constant $\psi$ is a solution of (\ref{eq:condpsi}). 

Saying that $u_f = u$ is equivalent to saying that $\int f (v-u) dv = 0$. So $\psi$ defined by (\ref{eq:condpsi}) is such that the following implication holds: 
\begin{equation}\label{E:CI1}
\int f(v_k-u_k) dv = 0, \quad \forall k \in \{1,\ldots, d\}  \quad \Longrightarrow \quad \int R(u;f) \, \psi \, dv = 0.
\end{equation}
Since both sides of the implication of \eqref{E:CI1} are linear forms of $f$, by a standard theorem \cite{Brezis}, there exists $\beta_u \in \mathbb{R}^d$ such that 
\begin{align*}
\int R(u;f) \, \psi \, dv 
&=  \beta_u \cdot \int f(v-u) dv.
\end{align*}
Using (\ref{E:collision_u}) and Green's formula (see (\ref{eq:self0})), and introducing the change of function $g = f /M_u$, it follows that $\psi$ is a solution of the following problem:
\begin{equation}
\int M_{u} \nabla_v g \nabla_v \psi \, dv = \int g \beta'_u \cdot (v-u) \, M_u \, dv, \quad \forall g \in V_u,
\label{eq:var}
\end{equation}
with $\beta'_u = - \beta_u / D$. We will drop the primes in the following.

Let $\phi = \beta_u \cdot (v-u)$. Then, problem (\ref{eq:var}) for $\psi$ can be equivalently written according to the variational formulation: 
\begin{equation} (\psi, g)_{V_u} = (\phi, g)_{H_u}, \quad \forall g \in V_u.
\label{E:VP} 
\end{equation}
Since any constant $\psi$ is a solution of (\ref{eq:condpsi}), we can subtract $\int \psi M_u dv$ from $\psi$ and assume, without loss of generality, that $\psi$ is such that $\int \psi \, M_u \, dv = 0$. Next, we have: 
\begin{equation}\int \phi \, M_u \, dv = \beta_u \cdot \int (v-u) \, M_u \, dv = 0, 
\label{eq:intphi}
\end{equation}
by the definition of $M_u$. Therefore, if $g$ is a constant, we have both $(\psi, g)_{V_u} = 0$ (by the definition (\ref{eq:inner}) of $(\cdot, \cdot)_{V_u}$) and $(\phi, g)_{H_u} = 0$ (by (\ref{eq:intphi})). Therefore, it is possible to restrict (\ref{E:VP}) to functions $g$ such that $(g,1)_{H_u}=0$. Consequently, we define the following space: 
$$ \dot{V}_u := \{ \phi \in V_u : \int \phi \, M_u \, dv = 0 \}, $$
and variational formulation (\ref{E:VP}) can be made precise as follows: 
\begin{equation} 
\mbox{Find } \psi \in \dot{V}_u \mbox{ such that } (\psi, g)_{V_u} = (\phi, g)_{H_u}, \quad \forall g \in \dot V_u.
\label{E:VP2} 
\end{equation}

The Poincar\'e inequality for Gaussian measures shows that for all $\phi \in V_u$,
\begin{equation}\label{E:grossIneq}
|\phi|_{V_u}^2 + (\phi,1)_{H_u}^2 \geq |\phi|_{H_u}^2. 
\end{equation}
Then following from (\ref{E:grossIneq}), the bilinear form $(\cdot, \cdot)_{V_u}$ is coercive on $\dot V_u$. By the Lax-Milgram theorem, we deduce that there exists a unique $\psi \in \dot V_u$ such that (\ref{E:VP2}) holds. But on the other hand, obvious calculation shows that $\psi(v) = T \beta_u \cdot (v-u)$ is a particular solution. Since it belongs to $\dot V_u$, it is the only solution of the variational formulation (\ref{E:VP2}). 

Adding any constant, we have just shown that any solution of (\ref{eq:condpsi}) is of the form $\psi(v) = \alpha + \beta \cdot v$, where $\alpha \in {\mathbb R}$ and $\beta \in {\mathbb R}^d$ are arbitrary. Since this form is independent of $u$, we deduce that such $\psi$'s are solutions of (\ref{eq:condpsi}) for all $u \in {\mathbb R}^d$ and are therefore the only CI's. This ends the proof of Proposition \ref{prop:CI}. \endproof

\section{Appendix: Expansion of the force term}
\label{app_expan}

{\bf Proof of Lemma \ref{lem:expanFf}.} By the change of variables $y = x - \varepsilon R \, \xi$ and Taylor's formula, we have 
\begin{eqnarray}
& & \hspace{-1cm}\int K \big(\frac{|x-y|}{\veps R} \big) \, \omega \, f(y,\omega,t)\,dy \,d\omega = \varepsilon^d \int K(|\xi|) \, \omega \,  f(x-\varepsilon \xi ,\omega,t)\,d\xi \,d\omega \nonumber\\
& & \hspace{-1cm} = (\varepsilon R)^d \int K(|\xi|) \, \omega \left[  f - \varepsilon R \, \nabla_x f \cdot \xi + \frac{\varepsilon^2 R^2}{2} D^2_x f : (\xi \otimes \xi) + O(\varepsilon^3) \right](x,\omega,t)
\,d\xi \,d\omega \nonumber \\
& & \hspace{-1cm} = (\varepsilon R)^d \, \big(j(x,t) + \varepsilon^2 R^2 k \, \Delta_x j + O(\varepsilon^4) \big). \label{eq:Kj}
\end{eqnarray}
We have used the definition of $k$ and the evenness of $K$ with respect to $\xi$ in order to cancel the odd order terms of the expansion. We have denoted by $D^2_x f$ the Hessian matrix of $f$ with respect to $x$ (i.e. the matrix of the second order derivatives) and the symbol `:' refers to the contracted product of tensors. In a similar way, we have:
\begin{eqnarray}
& & \hspace{-1cm}\int K \big(\frac{|x-y|}{\veps R}\big) f(y,\omega,t)\,dy \,d\omega =  (\varepsilon R)^d \big(\rho(x,t) + \varepsilon^2 R^2 k \,  \Delta_x \rho + O(\varepsilon^4) \big). \label{Krho}
\end{eqnarray}
Therefore, by expanding the ratio of (\ref{eq:Kj}) and (\ref{Krho}) up to the fourth order and in view of (\ref{vfeps}), we find (\ref{eq:vfeps2}) with (\ref{eq:uf1}). Formula (\ref{eq:Ffeps2}) immediately follows from (\ref{vfeps}). \endproof


\section{Appendix: proof of Theorem \ref{thm:NS}}
\label{app_NS_proof}

Multiplying (\ref{E:ourPDE_10}) by the collision invariants $1$ and $v$, we are led to the conservation equations (in the same way as in Section \ref{sec:hydro}): 
\begin{eqnarray} \label{E:massConservation_CE}
& & \partial_t \rho^\veps + \nabla_x \cdot j^\veps = 0, \\
\label{E:conservationOfMomentum_CE}
& & \partial_t j^\veps + \nabla_x \cdot \Sigma^\veps = \tau^{-1} q^\varepsilon + \varepsilon \sigma^{-1} r^\varepsilon, 
\end{eqnarray}
with $\rho^\veps$, $j^\veps = \rho^\varepsilon u^\varepsilon$, $\Sigma^\veps$, $q^\varepsilon$ given by (\ref{eq:rhojSeps}) and (\ref{eq:qeps}) and 
\begin{equation} r^\varepsilon = \int u_{f^\varepsilon}^1 f^\varepsilon \, dv = u_{f^\varepsilon}^1 \rho ^\varepsilon,
\label{eq:reps}
\end{equation}
(since $u_{f^\varepsilon}^1$ does not depend on $v$). 
We also note that we can write
\begin{equation}  \Sigma^\veps = \rho^\veps (u^\varepsilon \otimes u^\varepsilon) + S^\varepsilon, \quad 
S^\varepsilon = \int f^\varepsilon \, (v-u^\varepsilon) \otimes (v-u^\varepsilon) \, dv. 
\label{eq:Sigmaeps}
\end{equation}

The Chapman-Enskog expansion consists in closing the expressions of $S^\varepsilon$, $q^\varepsilon$ and $r^\varepsilon$ by a first order expansion of $f^\varepsilon$. To this end, we write the so-called macro-micro decomposition:
\begin{equation} 
f^\varepsilon = \rho^\varepsilon M_{u^\varepsilon} + \varepsilon f^\varepsilon_1.
\label{eq:macromicro}
\end{equation}
By the definition of $\rho^\varepsilon$ and $u^\varepsilon$, we have 
$$ \int f^\varepsilon \, dv = \int \rho^\varepsilon M_{u^\varepsilon} \, dv \quad \mbox{ and } \quad 
\int f^\varepsilon \, v \,  dv = \int \rho^\varepsilon M_{u^\varepsilon} \, v \, dv . $$
Consequently 
\begin{equation} 
\int f^\varepsilon_1 \, dv = 0, \quad \int f^\varepsilon_1 \, v \, dv = 0.
\label{eq:macromicro_2}
\end{equation}
The first term at the right-hand side of (\ref{eq:macromicro}) is the macroscopic part, as it is proportional to a local thermodynamical equilibrium and carries all information about the moments of the solution. The second term is the microscopic part. It carries no information about the macroscopic moments but instead carries information about the discrepancy between $f^\varepsilon$ and the local thermodynamical equilibrium. From Section \ref{sec:hydro}, we know that, provided that $\rho$ and $u$ satisfy the Euler equations, the microscopic part is small of order $\varepsilon$. This is why this microscopic part is multiplied by $\varepsilon$ in (\ref{eq:macromicro}). We stress the fact that there is no approximation involved (at this step) in (\ref{eq:macromicro}): it is a mere definition of $f^\varepsilon_1$. 

Inserting (\ref{eq:macromicro}) into the formulas providing the expressions of the quantities involved in the moment equations (\ref{E:massConservation_CE}), (\ref{E:conservationOfMomentum_CE}) (specifically equations (\ref{eq:Sigmaeps}), (\ref{eq:qeps}) and (\ref{eq:reps})), we can write the moment equations as follows: 
\begin{align}
&\partial_t \rho^\varepsilon + \nabla_x \cdot (\rho^\varepsilon u^\varepsilon ) = 0 , \label{E:alpha1_CE}\\
&\partial_t (\rho^\varepsilon u^\varepsilon ) + \nabla_x \cdot (\rho^\varepsilon u^\varepsilon \otimes u^\varepsilon ) + T \nabla_x \rho^\varepsilon = -\frac{1}{\tau} \rho^\varepsilon u^\varepsilon \big(\frac{|u^\varepsilon|^2 + (d+2) T}{a^2} - 1 \big) \nonumber \\
& \hspace{7cm} + \varepsilon ({\mathcal B}_1^\varepsilon + {\mathcal B}_2^\varepsilon + {\mathcal B}_3^\varepsilon) + O(\varepsilon^2) , \label{E:alpha2_CE}
\end{align}
with 
\begin{eqnarray}
& & {\mathcal B}_1^\varepsilon = \tau^{-1} \int \big( 1 - \frac{|v|^2}{a^2} \big) v f^\varepsilon_1 \, dv, \label{eq:B1} \\
& & {\mathcal B}_2^\varepsilon = \sigma^{-1} \rho^\varepsilon \, u^1_{\rho^\varepsilon M_{u^\varepsilon}} = \frac{k_R}{\sigma} (2 (\nabla_x \rho^\varepsilon \cdot \nabla_x) u^\varepsilon + \rho^\varepsilon \Delta_x u^\varepsilon), \label{eq:B2} \\
& & {\mathcal B}_3^\varepsilon = - \nabla_x \cdot \left( \int f_1^\varepsilon \, (v-u^\varepsilon) \otimes (v-u^\varepsilon) \, dv \right), \label{eq:B3}
\end{eqnarray}
where we have used Green's formula for ${\mathcal B}_1^\varepsilon$ and (\ref{eq:uf1}), (\ref{eq:reps}) for ${\mathcal B}_2^\varepsilon$. Now, in order to compute ${\mathcal B}_1^\varepsilon$ and ${\mathcal B}_3^\varepsilon$, we need to evaluate $f^\varepsilon_1$. But since we only look for ${\mathcal O}(\varepsilon)$ correction terms, and $f_1^\varepsilon$ is multiplied by $\varepsilon$, we may compute $f_1^\varepsilon$ up to terms of order ${\mathcal O}(\varepsilon)$.

Inserting (\ref{eq:macromicro}) into (\ref{E:collision2}), the collision operator $Q$ can be written:
\begin{equation} 
Q(f^\varepsilon) = \varepsilon D \nabla_v \cdot \left[ M_{u^\varepsilon} \nabla_v \big( \frac{f^\varepsilon_1}{M_{u^\varepsilon}} \big) \right] . 
\label{eq:Qnew}
\end{equation}
Inserting it in (\ref{E:ourPDE_10}), we get:
\begin{eqnarray} 
&&\hspace{-1cm}\partial_t f^\veps + \nabla_x \cdot (v f^\veps) + \tau^{-1} \nabla_v \cdot ( (1-\frac{|v|^2}{a^2})v f^\veps) + \varepsilon \sigma^{-1} \nabla_v \cdot (u_{f^\varepsilon}^1 f^\varepsilon) = \nonumber \\
&&\hspace{6.cm}
= D \nabla_v \cdot \left[ M_{u^\varepsilon} \nabla_v \big( \frac{f^\varepsilon_1}{M_{u^\varepsilon}} \big) \right], \label{E:ourPDE_11}
\end{eqnarray}
But we can neglect all terms of order $\varepsilon$ or more in (\ref{E:ourPDE_11}). Therefore, we are led to the following equation for $f^\varepsilon_1$: 
\begin{align} \label{E:ourPDE_12}
D \nabla_v \cdot \left[ M_{u^\varepsilon} \nabla_v \big( \frac{f^\varepsilon_1}{M_{u^\varepsilon}} \big) \right] = {\mathcal R}^\epsilon + O(\varepsilon).  
\end{align}
where
\begin{align} \label{E:ourPDE_13}
{\mathcal R}^\epsilon = \partial_t (\rho^\varepsilon M_{u^\varepsilon}) + \nabla_x \cdot (v \rho^\varepsilon M_{u^\varepsilon}) + \tau^{-1} \nabla_v \cdot ( (1-\frac{|v|^2}{a^2})v \rho^\varepsilon M_{u^\varepsilon}) .  
\end{align}
The inversion of (\ref{E:ourPDE_12}) will give us $f^\varepsilon_1$. Once $f^\varepsilon_1$ is obtained, we insert it into (\ref{eq:B1})-(\ref{eq:B2}) and this leads us to the expressions of the Navier-Stokes terms. We first compute the right-hand side ${\mathcal R}^\epsilon$: 

\begin{lemma}
We have (dropping the superscript $\varepsilon$ for the sake of clarity):
\begin{eqnarray}
& & \hspace{-1cm} {\mathcal R} = - \rho M_u \Big\{ h(v-u) :  \nabla_x u   + \frac{1}{\tau a^2} \big[ (d+2) T b(v-u) - a^2  d c(v-u)  \nonumber  \\
& & \hspace{0.5cm}  + (d+2)  e(v-u):(u \otimes u) + 3 T^{1/2} (d+2) g(v-u) \cdot u \big] \Big\} + {\mathcal O}(\varepsilon) ,
\label{eq:calR}
\end{eqnarray}
where $b(w)$, $c(w)$ are scalars, $e(w)$, $h(w)$ are tensors and $g(w)$ is a vector, and are given by:
\begin{eqnarray}
& & \hspace{-1.5cm} h(w) = \mbox{Id} - \frac{w \otimes w}{T} ,\quad b(w) = \frac{|w|^2}{T} \big( 1 - \frac{|w|^2}{(d+2)T} \big) , \quad c(w) =  1 - \frac{|w|^2}{dT} , \label{eq:hbc} \\
& & \hspace{-1.5cm} e(w) = \big( 1 - \frac{|w|^2}{(d+2)T} \big) \mbox{Id} - \frac{2 w \otimes w}{(d+2)T} , \quad g(w) = \big( 1 - \frac{|w|^2}{(d+2)T} \big) \, \frac{w}{T^{1/2}} .
\label{eq:eg} 
\end{eqnarray}
By construction, these quantities are dimensionless. 
\label{lem;Reps}
\end{lemma}

\medskip
\noindent
{\bf Proof.} We use the hydrodynamic equations (\ref{E:alpha1_CE}), (\ref{E:alpha2_CE}) (dropping the $O(\varepsilon)$ terms), in order to replace the time derivatives by space derivatives in (\ref{E:ourPDE_13}). This procedure is a classical step of any Chapman-Enskog expansion. For simplicity of notation, we omit the dependencies on $\varepsilon$. Concerning the first two terms of (\ref{E:ourPDE_13}), we write:
\begin{eqnarray*}
& & \hspace{-0.5cm}(\partial_t + v \cdot \nabla_x) (\rho M_u) = M_u \left\{ (\partial_t + v \cdot \nabla_x) \rho + \rho  (\partial_t + v \cdot \nabla_x) (\ln M_u) \right\} \\
& &  = M_u \left\{ (\partial_t + v \cdot \nabla_x) \rho + \rho \frac{v-u}{T} \cdot  (\partial_t + v \cdot \nabla_x) u \right\} \\
& &  = M_u \left\{ \partial_t \rho + u \cdot \nabla_x \rho + (v-u) \cdot \nabla_x \rho  \right. \\
& & \hspace{3cm} \left. + \rho \frac{v-u}{T} \cdot  (\partial_t u + (u \cdot \nabla_x) u + ((v- u) \cdot \nabla_x) u )\right\}  .
\end{eqnarray*}
Now, using (\ref{E:alpha1_CE}), (\ref{E:alpha2_CE}) (dropping the $O(\varepsilon)$ terms), we have: 
\begin{eqnarray*}
& & \partial_t \rho + u \cdot \nabla_x \rho = - \rho \nabla_x \cdot u + O(\varepsilon) , \\
& & \partial_t u + u \cdot \nabla_x u = - T \nabla_x \ln \rho - \tau^{-1} u \big(\frac{|u|^2 + (d+2) T}{a^2} - 1 \big) + O(\varepsilon) . 
\end{eqnarray*}
Inserting these expressions into the previous ones leads to: 
\begin{eqnarray}
& & \hspace{-1cm}(\partial_t + v \cdot \nabla_x) (\rho M_u) = M_u \Big\{ - \rho \nabla_x \cdot u  + (v-u) \cdot \nabla_x \rho  \nonumber\\
& & \hspace{1.5cm} + \rho \frac{v-u}{T} \cdot  \big[ - T \nabla_x \ln \rho - \tau^{-1} u \big(\frac{|u|^2 + (d+2) T}{a^2} - 1 \big)  \nonumber \\
& & \hspace{6.1cm} +
((v- u) \cdot \nabla_x) u \big] \Big\} + O(\varepsilon) \nonumber\\
& & \hspace{-1cm} = M_u \Big\{ - \rho \nabla_x \cdot u  + \rho \frac{v-u}{T} \cdot  \big[ - \tau^{-1} u \big(\frac{|u|^2 + (d+2) T}{a^2} - 1 \big)  \nonumber \\
& & \hspace{6.1cm} + ((v- u) \cdot \nabla_x) u \big] \Big\} + O(\varepsilon) \nonumber\\
& & \hspace{-1cm} = M_u \Big\{ - \rho h(v-u) :  \nabla_x u  - \frac{\rho}{\tau T} \big(\frac{|u|^2 + (d+2) T}{a^2} - 1 \big)  (v-u) \cdot u  \Big\} + O(\varepsilon) \nonumber\\
\label{terms12}
\end{eqnarray}

For the third term (\ref{E:ourPDE_13}), we compute
\begin{eqnarray}
&& \hspace{-1.cm} \tau^{-1} \nabla_v \cdot ( (1-\frac{|v|^2}{a^2})v \rho M_{u}) = \nonumber \\
&& \hspace{1.5cm} - \tau^{-1} \rho M_u \left\{ (d+2) \frac{|v|^2}{a^2} + \big( 1 - \frac{|v|^2}{a^2} \big) \frac{v \cdot (v-u)}{T} - d \right\} .
\label{eq:third}
\end{eqnarray}
Collecting the second term at the right-hand side of (\ref{terms12}) together with (\ref{eq:third}) leads to an expression ${\mathcal S}$ which we can split in the following way: 
\begin{eqnarray}
& & \hspace{-0.5cm} {\mathcal S} = - M_u  \frac{\rho}{\tau T} \big(\frac{|u|^2 + (d+2) T}{a^2} -1 \big)  (v-u) \cdot u  + \frac{1}{\tau} \nabla_v \cdot ( (1-\frac{|v|^2}{a^2})v \rho M_{u}) \nonumber \\
& & = - \frac{\rho}{\tau a^2} M_u \left\{ (d+2) T b(v-u) - a^2 d c(v-u) \right. \nonumber  \\
& & \hspace{1.5cm} \left. + (d+2) e(v-u):(u \otimes u) + 3 T^{1/2} (d+2) g(v-u) \cdot u \right\} ,
\label{eq:calS}
\end{eqnarray}
Collecting this expression with the first term at the right-hand side of (\ref{terms12}) leads to (\ref{eq:calR}). \endproof

In order to solve equation (\ref{E:ourPDE_12}) for $f^\varepsilon_1$, we need to solve equations of the type 
\begin{equation}
L_uf := - D \nabla_v \cdot \big[ M_u \nabla_v \big( \frac{f}{M_u} \big) \big] = g, \label{eq:elliptic}
\end{equation}
where $u$ is an arbitrary vector of ${\mathbb R}^d$ and $g$ is a given function. We refer the reader to Section \ref{subsec:CI} for the definitions of the spaces $H_u$, $V_u$ and $\dot V_u$. We state the following lemma, whose proof is identical to that of proposition \ref{prop:CI} and is left to the reader. 

\begin{lemma}
Let $g$ be such that $g/M_u \in H_u$. Then, equation (\ref{eq:elliptic}) has a solution if and only if $g$ satisfies the solvability condition
\begin{equation}
\int g \, dv = 0 .  \label{eq:solv_g}
\end{equation}
Under this condition, problem (\ref{eq:elliptic}) has a unique solution $f$ such that $f/M_u \in \dot V_u$, or in other words, such that $f/M_u \in V_u$ and satisfies 
\begin{equation}
\int f \, dv = 0 .  \label{eq:solv_f}
\end{equation}
This unique solution $f$ is denoted by $f = L_u^{-1} g$ and $L_u^{-1}$ is called the pseudo-inverse of $L_u$. The set of solutions to equation (\ref{eq:elliptic}) is given by $L_u^{-1} g + \mbox{Span}(M_u)$. Additionally, if $g$ is such that 
\begin{equation}
\int g \, v \, dv = 0 ,  \label{eq:solv_gv}
\end{equation}
then, $f = L_u^{-1} g$ satisfies 
\begin{equation}
\int f \, v \, dv = 0 .  \label{eq:solv_fv}
\end{equation}
\label{lem:inversion_Lu}
\end{lemma}

\medskip
\noindent
Now, we verify that each of the elementary functions $b(v-u)M_u$, $c(v-u)M_u$, $e(v-u)M_u$, $h(v-u)M_u$ and $g(v-u)M_u$ satisfy both (\ref{eq:solv_g}) and (\ref{eq:solv_gv}) and therefore, that the corresponding equation (\ref{eq:elliptic}) is invertible. More precisely, we have:

\begin{lemma}
The functions $h(v-u)M_u$, $b(v-u)M_u$, $c(v-u)M_u$, $e(v-u)M_u$  and $g(v-u)M_u$ satisfy  (\ref{eq:solv_g}) and (\ref{eq:solv_gv}). We introduce $ H(v-u) M_u = -L_u^{-1} (h(v-u)M_u)$ and similarly for $B$, $C$, $E$ and $G$. We have: 
\begin{eqnarray}
& & H(w) = - \frac{\sigma}{2} h(w) , \quad B(w) = - \frac{\sigma d}{4} (c(w) + \frac{1}{d} b(w)) , \label{eq:HB} \\
& & C(w) =  - \frac{\sigma}{2} c(w) , \quad  E(w) = - \frac{\sigma}{2} e(w) , \quad G(w) =  - \frac{\sigma}{3} g(w),  \label{eq:CEG} 
\end{eqnarray}
and $ H(v-u) M_u$ through $G(v-u) M_u$ satisfy  (\ref{eq:solv_f}) and (\ref{eq:solv_fv}). Then: 
\begin{eqnarray}
& & \hspace{-1cm} f_1^\varepsilon = \sigma \rho M_u \Big\{ \frac{1}{2} h(v-u) :  \nabla_x u   \nonumber \\
& & \hspace{0.5cm} + \frac{1}{\tau a^2} \Big[ \frac{(d+2) T}{4} b(v-u) + \frac{d}{2} \big( \frac{(d+2) T }{2} - a^2 \big)  c(v-u)  \nonumber  \\
& & \hspace{0.5cm}  + \frac{(d+2)}{2}  e(v-u):(u \otimes u) + T^{1/2} (d+2) g(v-u) \cdot u \Big] \Big\} + {\mathcal O}(\varepsilon)  .
\label{eq:f10_2}
\end{eqnarray}
\label{lem:f10}
\end{lemma}

\medskip
\noindent
{\bf Proof.} The proof that $h(v-u)M_u$ through $g(v-u)M_u$ satisfy  (\ref{eq:solv_g}) and (\ref{eq:solv_gv}) easily follows from classical formulas for moments of the Gaussian, which we leave to the reader. Then, we apply Lemma \ref{lem:inversion_Lu}, which gives the existence of $L_u^{-1} (h(v-u)M_u)$ through $L_u^{-1} (g(v-u)M_u)$ and the fact that they satisfy (\ref{eq:solv_f}) and (\ref{eq:solv_fv}). Formulas (\ref{eq:HB}), (\ref{eq:CEG}) follow from explicitly computing the action of $L_u^{-1}$ and using the uniqueness statement of Lemma \ref{lem:inversion_Lu}. Finally, equation (\ref{E:ourPDE_12}) for $f_1^\varepsilon$, which can be written (up to order $O(\varepsilon)$ terms) $-L_u f_1^\varepsilon = {\mathcal R}$, can be solved by $f_1^\varepsilon = - L_u^{-1} {\mathcal R}$ since according to the first equation of (\ref{eq:macromicro_2}), $f_1^\varepsilon$ satisfies (\ref{eq:solv_f}). By the linearity of $L_u$ and the decomposition (\ref{eq:calR}) of ${\mathcal R}$, we can write:
\begin{eqnarray}
& & \hspace{-1cm} f_1^\varepsilon= - \rho M_u \Big\{ H(v-u) :  \nabla_x u   + \frac{1}{\tau a^2} \big[ (d+2) T B(v-u) - a^2  d C(v-u)  \nonumber  \\
& & \hspace{0.5cm}  + (d+2)  E(v-u):(u \otimes u) + 3 T^{1/2} (d+2) G(v-u) \cdot u \big] \Big\} + {\mathcal O}(\varepsilon) .
\label{eq:f_100}
\end{eqnarray}
Thanks to (\ref{eq:HB}), (\ref{eq:CEG}), equation (\ref{eq:f10_2}) follows. \endproof

\medskip
\noindent
We are now in a position to calculate ${\mathcal B}_1$ and ${\mathcal B}_3$ (see (\ref{eq:B1}), (\ref{eq:B3})). We state: 

\begin{lemma}
We have:  
\begin{eqnarray}
& & \hspace{-1cm} {\mathcal B}_1 = \frac{\lambda}{2} \rho \Big\{ \big[ (\nabla_x \cdot u) u + \nabla_x \big( \frac{|u|^2}{2} \big) + (u \cdot \nabla_x) u \big] + \frac{d+8}{\tau} u \big[\frac{|u|^2}{a^2}  - \nu \big] \Big\} + {\mathcal O}(\varepsilon), \label{express_B1} \\
& & \hspace{-1cm} {\mathcal B}_3 = \mu \nabla_x  \cdot (\rho {\mathcal E}(u)) +  \nabla_x \pi(\rho,u) + \lambda \nabla_x \cdot (\rho u \otimes u) + {\mathcal O}(\varepsilon), \label{express_B3} 
\end{eqnarray}
with ${\mathcal E}(u)$, $\mu$, $\lambda$ and $\pi(\rho, u)$ are defined at theorem \ref{thm:NS} and $\nu$ is given by: 
\begin{eqnarray}
& &  \nu = \frac{d+2}{d+8}\big( 1 - (d+4) \frac{T}{a^2} \big) ,  \label{eq:nu} 
\end{eqnarray}
\label{lem:B1-B3}
\end{lemma}

\medskip
\noindent
{\bf Proof.} We first consider ${\mathcal B}_1$. Splitting $v$ into $(v-u) = u$ and using (\ref{eq:macromicro_2}), we have
\begin{eqnarray*}
{\mathcal B}_1 &=& - \frac{1}{\tau a^2} \int \big\{ |v-u|^2 (v-u) + |v-u|^2 u + 2 ((v-u) \otimes (v-u)) u\big\}  f_1^\varepsilon \, dv \\
&=& J_1 + J_2 + J_3.  
\end{eqnarray*}
$J_1$ involves an integral of an odd power of $v-u$ and only the $g$ term in $ f_1^\varepsilon$ contributes to it. $J_2$ and $J_3$ involve integrals of even powers of $v-u$ and therefore, only the $h$, $b$, $c$ and $e$ terms need to be taken into account. The computation of these terms rely on computing moments of the Gaussian which are left to the reader. We find:
\begin{eqnarray*}
J_1 &=& \frac{2 T^2 \sigma (d+2)}{\tau^2 a^4} \rho u + {\mathcal O}(\varepsilon), \\
J_2 &=& \frac{T \sigma}{\tau a^2}  (\nabla_x \cdot u) \rho u + \frac{d \sigma T}{\tau^2 a^2} \big[\frac{(d+2)T}{a^2} -1 \big] \rho u + \frac{(d+2) \sigma T}{\tau^2 a^4} \rho |u|^2 u + {\mathcal O}(\varepsilon),\\
J_3 &=& \frac{T \sigma}{\tau a^2} \rho ((\nabla_x u) u + (\nabla_x u)^T u) +  \frac{2 \sigma T}{\tau^2 a^2} \big[\frac{(d+2)T}{a^2} -1 \big] \rho u  \\
& & \hspace{7cm} + \frac{6 \sigma T}{\tau^2 a^4} \rho |u|^2 u + {\mathcal O}(\varepsilon).
\end{eqnarray*}
Now, adding up these three expressions leads to (\ref{express_B1}).

We now turn our attention to ${\mathcal B}_3$. For this purpose, we compute the tensor 
\begin{eqnarray*}
{\mathcal U} = \int f_1^\varepsilon \, (v-u^\varepsilon) \otimes (v-u^\varepsilon) \, dv. 
\end{eqnarray*}
Since the integral involves an even power of $v-u$, only the $h$, $b$, $c$ and $e$ terms of $f_1^\varepsilon$ need to be taken into account. The computation leads to 
\begin{eqnarray*}
&& \hspace{-1cm} {\mathcal U} = - \frac{\sigma T}{2} \rho ((\nabla_x u) + (\nabla_x u)^T) - \rho \frac{\sigma T}{\tau} \big[\frac{(d+2)T}{a^2} -1 \big] \mbox{Id} \\
&& \hspace{5cm}  - \frac{\sigma T}{\tau a^2} \rho (|u|^2 \mbox{Id} + 2 u \otimes u) + {\mathcal O}(\varepsilon). 
\end{eqnarray*}
Inserting this expression into (\ref{eq:B3}) leads to (\ref{express_B3}). \endproof

Now, by adding the expressions of ${\mathcal B_1}$ through ${\mathcal B_3}$ found above into (\ref{E:alpha2_CE}), we find the following momentum equation:   
\begin{align}
&\partial_t (\rho u ) + \nabla_x \cdot (\rho u \otimes u ) + T \nabla_x \rho = -\frac{1}{\tau} \rho u \big(\frac{|u|^2 + (d+2) T}{a^2} - 1 \big) \nonumber \\
& \hspace{1cm} + \varepsilon \frac{\lambda}{2} \rho \Big\{ \big[ (\nabla_x \cdot u) u + \nabla_x \big( \frac{|u|^2}{2} \big) + (u \cdot \nabla_x) u \big] + \frac{d+8}{\tau} u \big[\frac{|u|^2}{a^2}  - \nu \big] \Big\} \nonumber \\
& \hspace{1cm} + \varepsilon \frac{k_r}{\sigma} (2 (\nabla_x \rho \cdot \nabla_x) u + \rho \Delta_x u) \nonumber \\
& \hspace{1cm} + \varepsilon \mu \nabla_x  \cdot (\rho {\mathcal E}(u)) +  \nabla_x \pi(\rho,u) + \lambda \nabla_x \cdot (\rho u \otimes u) + {\mathcal O}(\varepsilon^2).  \label{E:alpha2_CE12}
\end{align} 
Now, the second term at the left-hand side of (\ref{E:alpha2_NS}) combines the second term at the left-hand side and the last term of the right-hand side of (\ref{E:alpha2_CE12}); the third term at the left-hand side of (\ref{E:alpha2_NS}) combines the third term at the left-hand side and the penultimate term of the right-hand side of (\ref{E:alpha2_CE12}); the first term at the right-hand side of (\ref{E:alpha2_NS}) combines the first term at the right-hand side and the last term of the second line of (\ref{E:alpha2_CE12}); and the other terms are unchanged but merely re-ordered. This ends the proof of Theorem \ref{thm:NS}. \endproof


\end{document}